# SOLUTE TRANSPORT IN A POROUS MEDIUM:
# A MASS-CONSERVING SOLUTION FOR THE CONVECTION-DISPERSION EQUATION IN A FINITE DOMAIN

WILLIAM J. GOLZ
PH.D. DISSERTATION, LOUISIANA STATE UNIVERSITY, 2003

**ACKNOWLEDGMENTS**

I owe a special note of thanks to my major advisor from Civil and Environmental Engineering, Dr. D. Dean Adrian, for introducing me to the broad literature of solute transport, and to my minor advisor from Mathematics, Dr. J. Robert Dorroh, for encouraging me to investigate the subtleties in my mathematical statements.

I appreciate the thoughtful evaluation of my dissertation provided by the Dean's representative, Dr. Gary A. Breitenbeck from Agronomy, and by my committee members from Civil and Environmental Engineering, Dr. Kelly Ann Rusch and Dr. Vijay P. Singh.

My greatest debt for the fruition of this dissertation, as for most good things in my life, is to my wife Annette, for her faith.

**TABLE OF CONTENTS**

**LIST OF FIGURES**

**ABSTRACT**

This dissertation considers the proper mathematical description for the physical problem of a miscible solute undergoing longitudinal convective-dispersive transport with constant production, first-order decay, and equilibrium sorption in a porous medium. Initial and input concentrations may be any continuously differentiable functions and the mathematical system is articulated for a finite domain. This domain yields a mass balance which requires Robin (i.e., third-type) boundaries, which describe a continuous flux but a discontinuous resident-concentration. The discontinuity in the resident concentration at the outflow boundary yields an underdetermined system when the exit concentration is not experimentally measured. This is resolved by defining the unknown effluent concentration from a semi-infinite problem which satisfies a Dirichlet (i.e., first-type) condition at the origin.

The solution is represented in a uniformly convergent series of real variables. The representation can be sequenced to describe any configuration of discrete reactors or approach reservoirs. Individual reacting segments are allowed to have differing lengths and transport parameters up to the complexity of the governing equation. Such discrete segments may be constructed from finitely small slices to approximate a continuous variation in any of the modeled parameters, such as velocity or diffusion. The physical phenomenon that can be described include layered hydrogeologic strata, as well as two- or three- dimensional transport when hydrodynamic properties exhibit a spatial proportionality.

The large volume of antecedent literature on finite solutions for convective-dispersive transport equations grew out of the historical precedents set by Danckwerts (1953) and Wehner and Wilhelm (1956) whom made simplifying assumptions of continuous boundary concentrations. This dissertation includes the demonstration that continuous-concentration hypotheses, whether rendered as Dirichlet or homogeneous Neumann (i.e., second-type) conditions, satisfy external mass conservation yet fail to provide solutions that are internally consistent with the governing equation.

## CHAPTER 1. INTRODUCTION

### 1.1. BRIEF HISTORY OF THE PROBLEM

Over the latter half of the past century our understanding of chemical processes has advanced by virtue of a steady stream of research that has defined much of the accepted theory that governs chemical-transformations, both pure and biologically mediated. While research at the beginning of this period was free to focus on industrial-chemical production, our science is now required to encompass an understanding of the behavior of dissolved contaminants in natural systems. A problem of continuing practical importance, common to the production and remediation of chemicals and therefore a focal point in the interdisciplinary literature, is that of a solute undergoing one-dimensional, convective-dispersive transport in a porous medium.

Physical investigations involving porous media present a unique challenge because solute behavior can change dramatically at a physical boundary. In addition, it is often impossible to measure interior concentrations without disturbing natural flow lines and thereby introducing error. These inherent complications mean that empirical observations will typically furnish only input and exit concentrations. Consequently, the complex physics that control solute behavior within the medium can be understood only through carefully defined mathematical idealizations. Those models then integrate, and become an integral part of, our scientific knowledge. If they are to furnish the foundation for a scientific dialogue, whether quantitative or qualitative, we must ask that our models remain true to accepted physical theory and mathematical logic.

Problems describing non-steady chemical transport in a finite domain lead naturally to a mass balance that requires a third-type condition at each boundary. However, a departure from that conservation requirement grew out of a strong precedent provided by three papers published around the middle of the last century. In 1944, Hulburt summarily assumed a first-type upper boundary and postulated that, because any reaction should have gone to completion before the solute reached the exit, the condition there should be a zero gradient. Danckwerts' (1953) well-known paper accepted a third-type entrance but presented arguments similar to Hulburt's in support of the zero-gradient exit. Another widely prescribed set of assumptions originated with Wehner and Wilhelm's (1956) work which admitted dispersion in the appended reservoirs by specifying a continuous concentration at both boundaries. These historic works concerned steady-state reactions, yet the specification of continuous concentrations has become customary for non-steady reactions (Parulekar and Ramkrishna 1984), and the zero-gradient exit is often accepted without comment (Liu et al. 2000).

### 1.2. ARTICULATION OF THE MATHEMATICAL SYSTEM

This dissertation proposes a proper mathematical description for the physical problem of a miscible solute undergoing longitudinal convective-dispersive transport with constant production, first-order decay, and equilibrium sorption in a porous medium. Our differential system will be derived from the demonstration that conservation of mass leads directly to the convection-dispersion equation (CDE)

$$RC_t = DC_{xx} - vC_x - \mu C + \gamma, \quad 0 \leq x \leq \ell, \quad t \in \mathbb{R}, \qquad (1.1)$$

and the Robin boundaries

$$vC(0,t) - DC_x(0,t) = vg(t), \quad t \geq t_0, \qquad (1.2)$$

$$vC(\ell,t) - DC_x(\ell,t) = vC_E, \quad t \geq t_0. \qquad (1.3)$$

We allow the general initial condition

$$C(x,t_0) = \phi(x), \quad 0 \leq x \leq \ell. \qquad (1.4)$$

The dependent variable $C = C(x,t)$ is a solute concentration which depends continuously upon distance (x) and time (t). The coefficients are constants which describe a retardation factor (R), hydrodynamic dispersion (D), and average interstitial fluid velocity (v), with production ($\gamma$) and decay ($\mu$) occurring in both the liquid and sorbed phase.[1] We ask that the initial concentration $\phi(x)$ be twice continuously differentiable and that time be restricted to an interval in which the solute input concentration $g(t)$ will be continuously differentiable.[2]

The exit concentration $C_E$ will be a measured quantity when, for example, the dispersion constant is to be derived from a time series of experimental measurements. However, even when the effluent concentration will not be defined from an empirical observation, the smoothing effect of the diffusive process yields the reasonable assumption that $C_E$ will be continuously differentiable. We thus regard $C_E$ as arbitrary and go on to transform (1.1)-(1.4) to an inhomogeneous diffusion equation with homogeneous boundaries. This yields a problem whose spatial operator defines a Sturm-Liouville Robin system which admits a series expansion in real variables that provides a solution to the original problem.

There are many problems of practical interest in which it is necessary to predict the exit concentration, for instance when designing experimental reactors or when planning sampling times to minimize the number of measurements required to accurately reconstruct dynamic effluent concentrations. In such cases, it will be necessary to have a formula for the exit concentration which we can furnish from the flux concentration

$$vC_F(x,t) = vC(x,t) - DC_x(x,t), \quad 0 \leq x < \infty. \qquad (1.5)$$

---

[1] For the decomposition of the coefficients, see sec. 3.2 especially eq. (3.4). For a list of the symbols and their corresponding units, see app. A.

[2] Although piecewise-continuous differentiability would be sufficient for our solution, continuous differentiability has the physical rationale that non-smooth initial concentrations typically exist only immediately after solute deposition or in the absence of an available solvent. Then, by considering time intervals in which the input concentration is continuously differentiable, we achieve the best possible convergence, and this has the practical advantage of admitting a strong justification for the analytic techniques, which translates into better numerical approximations.

Note that a careful substitution of (1.5) into (1.1) will yield the fact that $C_F(x,t)$ satisfies an equation identical to the CDE and that, by substituting (1.5) into (1.2),(1.4), we define the auxiliary conditions

$$v\,C_F(x,t_0) = v\,\phi(x) - D\,\phi_x(x), \qquad 0 \le x < \infty, \qquad (1.6)$$

$$C_F(0,t) = g(t), \qquad\qquad t \ge t_0. \qquad (1.7)$$

If we now refer to (1.3) for the fact that $C_F(\ell,t) = C_E$, it is clear that a solution for the flux concentration will also yield an expression for the exit concentration.

1.3. OBJECTIVE OF THIS DISSERTATION

The objective of this dissertation is to provide the proper solution to the mathematical system articulated in section 1.2. This will devolve upon two *a posteriori* demonstrations: (i) that the problem (1.1)-(1.4) is well posed, which means to show that we have constructed a representation for a solution which is unique and stable; and (ii) that the Robin boundaries are a necessary condition of any solution to the CDE.

**CHAPTER 2. LITERATURE REVIEW**

Mathematical solutions of the one-dimensional CDE provide descriptions for problems involving a variety of solute types and porous media (for a representative list, with sources, see van Genuchten and Alves 1982, p. 3). The expansive body of literature motivated by the broad applicability of the CDE has provided definitive answers to a number of questions concerning the physical phenomena described by different combinations of boundary conditions and governing equations for infinite domains (e.g., see Kreft and Zuber 1978). However, the customary boundary conditions for a finite domain, and the validity of all consequent solutions, have become an established source of profound and unresolved disagreement (see Parker 1984, Leij and Toride 1998). This review provides a synopsis of the central arguments in the deliberation over finite problems by discussing several of the most widely cited papers.

In 1944, Hulburt provided one of the earliest adaptations of hydrodynamic equations to chemical reactions in flowing systems. His primary objective was to relate steady-state chemical yields to equipment sizes and reaction rates. The formulas he developed described a homogeneous fluid subject to axial convection and dispersion with disappearance in a finite cylindrical vessel. At the entrance, a dimensionless unit concentration was permitted to enter the reactor unchanged according to a first-type condition. For the exit, he derived a third-type condition from the conservation of mass, but then argued that, because the reaction must cease at the outlet, this would reduce to a zero gradient (i.e., a homogeneous second-type condition).

In his 1953 paper, which provided the first complete treatment of residence times in flowing systems, Danckwerts discussed steady-state flow, dispersion, and decay in a packed bed. The discussion on page 10 of his paper acknowledged that mass conservation for a finite vessel led to a third-type condition for each boundary, and he accepted that condition for the entry. In his discussion of the exit, Danckwerts correctly dismissed a positive exit gradient with the argument that it would require the reactant to pass through a minimum in the interior of the vessel.[1] However, he goes on to reject the negative gradient on the grounds that the concentration just short of the exit plane would then be smaller than that in the exit stream; but that is precisely the discontinuity required to conserve mass. Danckwerts (1953) settled on a third-type entrance [ibid., eq. (31)] with a zero-gradient exit [ibid., eq. (32)]. Because his work is more widely recognized than any other in the topic area, his boundary specifications are commonly referred to as the Danckwerts' conditions.[2]

---

[1]The dismissal of a positive exit gradient is correct for Danckwerts' special case because there is no production within the reactor. However, for more general problems like the one we consider in this dissertation, the possibility of a positive gradient at the exit plane is entirely plausible and conforms to applicable maximum-minimum principles. For a discussion of those principles, cf. (i) Strauss (1992, sec. 2.3), or (ii) Gustafson (1999, problem 1.9.3).

[2]Parulekar and Ramkrishna (1984, p. 1571, footnote) cite a source which they claim traces the Danckwerts' boundary conditions to earlier work by Langmuir.

In 1956, Wehner and Wilhelm (WW) described a steady-state problem involving axial convection and dispersion with decay in a finite-length reactor with infinite fore and aft reservoirs. The distinguishing characteristic of their model was that it admitted dispersive reservoirs by equating third-type conditions on either side of the reactor boundaries for any nonzero, but possibly distinct, set of dispersion constants.[3]

For the upper end of their system, WW specified a normalized concentration of unity located at some large negative distance with a continuous concentration on either side of the inlet plane. Since the flux concentration is by definition continuous, their entry condition can be regarded as a specification of complete continuity. Hiby (1962) conducted a careful study of dispersion near an infusion point and found that the upstream mass-transport required to satisfy concentration continuity was inconsistent with experiment. Deckwer and Mählmann (1976) went on to study solute behavior near a sharp boundary to evaluate model assumptions. The authors trisected a reactor with sharp internal interfaces created by, for instance, different densities in adjoining sections. Their observations clearly confirmed that solute can exhibit a discontinuity when it crosses a plane where conditions change abruptly.[4]

WW noted that, within the reactor, their solution coincided with Danckwerts, and they went on to note that when dispersion was absent in the aft reservoir the limit of their exit condition was Danckwerts' zero gradient (WW, last paragraph of p. 91, and appendix). However, that zero gradient followed from the choice of problem geometry and boundary assumptions: the infinite aft reservoir meant that the exit concentration had to be nonincreasing so that the solution would remain bounded at a large distance [ibid., eq. (18), and discussion on p. 91]; then, the hypothesis that concentrations were continuous on either side of the exit plane [ibid., eq. (14)] carried the reservoir concentration into the reacting section.

In 1962, Bischoff and Levenspiel (BL) compared dispersion estimates from different models for a variety of solute-measurement points. The experimental conditions that they investigated involved non-steady axial convection with dispersion in a finite-length reactor with infinite entrance and exit reservoirs. For their principal model they adopted the conditions of WW with one exception: for the upper end of their system they specified a varying concentration located at some given distance above a fixed measurement point external to the entrance. BL developed expressions based upon lower measurement points located within the reactor, at the outflow boundary, and in the discharge reservoir. When the Danckwerts solution, which forced the exit gradient to zero, was compared to a semi-infinite solution, which was unaffected by the boundary, BL found that predicted Peclet numbers could differ by more than thirty percent.[5] Because BL made the *a priori* assumption that the finite exit conditions of WW and Danckwerts provided proper descriptions of boundary layer effects, the large discrepancies in Peclet numbers were attributed to an error in the semi-infinite expression.

---

[3]Standart (1968, p. 653) attributes the WW boundary conditions to earlier work by Damköhler.

[4]Deckwer and Mählmann (1976) summarized their results in their abstract by stating that a Danckwerts' solution was closer to experimental observation than a WW solution, but we should add that their results can only be regarded as a confirmation of Danckwerts third-type boundary. That is, the zero gradient is equivalent to the statement that concentrations are equal on either side of the exit plane, so the discontinuity is a contradiction.

[5]It is customary to define the Peclet number for this type of problem as the product of fluid velocity and length of flow domain divided by the dispersion constant, i.e., $v\ell/D$.

In 1984, van Genuchten and Parker (VGP) illustrated some important relationships between predicted concentrations and the spatial domain of a solution. They investigated a non-steady-state problem where convection, dispersion, and sorption occurred at small values of the column Peclet number. Their paper begins by comparing two solutions posed for a positive infinite domain that differed only in their entrance conditions. They demonstrated that a solution by Lindstrom et al. (1967) employing a third-type entrance would satisfy the mass-balance expression for a resident concentration (i.e., the concentration in resting fluid). It was then shown that a solution by Lapidus and Amundson (1952) which incorporated a first-type condition would describe a flux concentration (i.e., the rate of mass-flow per unit area). Those two solutions describe, respectively, the solute profile within a finite vessel and the inlet and exit concentrations just external to the boundary. This follows from the fact that the difference between resident and flux concentrations is proportional to the spatial derivative of the resident concentration (VGP, eq. [27], which transforms a resident concentration into one for flux).

VGP extended their discussion to two solutions for a finite domain which, like the infinite solutions, differed only in their entrance conditions. A solution by Brenner (1962), which employed a third-type entrance and a zero gradient at the exit (i.e., Danckwerts' conditions), was compared to a solution by Cleary and Adrian (1973), which employed a first-type entrance and a zero-gradient exit (i.e., Hulburt's conditions). VGP argued that, if the finite exit was correct, Cleary and Adrian's use of the first-type entrance should yield a flux concentration like the solution of Lapidus and Amundson. They went on to reason that when Brenner's expression employing a third-type entrance was transformed into an expression for flux it should have coincided with the solution of Cleary and Adrian throughout the domain. VGP then plotted a graph to demonstrate that the solutions agreed at the entrance but increasingly diverged as they approached the lower boundary (VGP, fig. 2).

Graphical analyses have found wide use in the literature as a tool to evaluate dissimilar expressions, but a graph can be misleading when there is hidden error, such as when solutions with different rates of convergence are treated similarly. For example, VGP superimpose the predicted-concentration curve of Brenner onto that of Lapidus and Amundson (VGP, fig. 3). Those curves are shown to differ substantially at a Peclet value of unity, and VGP regard this as further reason to dismiss Brenner's solution. However, this contravenes VGP's earlier demonstration that Brenner's solution satisfied the overall mass-balance for a finite column and must therefore provide correct exit concentrations (see VGP, eq. [22], and the following discussion).

The literature on finite reactors spans the latter half of the last century, and the first questions regarding the boundary conditions were raised soon after the earliest publications. These questions have remained unresolved, not only because of the strong historical precedents but because the comparative analysis of solutions cannot definitively answer the question of what constitutes a necessary condition for a solution to the CDE.

## CHAPTER 3. FORMULATION OF THE GOVERNING EQUATION AND BOUNDARY CONDITIONS

### 3.1. NOTE ON MATERIAL CONTINUA

A liquid continuum for a miscible solute is typically defined in terms of a representative elementary volume (REV) composed of a sphere with a diameter greater than the mean free path between molecules and sufficiently large to ensure that the behavior of any individual molecule is insignificant (see Bear 1972, sec. 1.3). This allows most molecular-scale behavior such as hydrodynamic dispersion and transformations to be described by constants which represent macroscopic averages so that the liquid can be regarded as homogeneous.

When a homogeneous liquid undergoes Darcy flow through a porous medium, the solid continuum is often defined in terms of a REV constituted by a sphere with a diameter sufficient to ensure that, wherever the sphere is placed, the volumetric porosity will remain constant. As demonstrated by Bear (1972, sec. 1.3.3), this leads to a single equivalent value for volumetric, areal, and linear porosity. Average porosity is normally a good index of all the physical characteristics of a solid that significantly affect hydrodynamics, e.g., particle size and shape.

For most practical problems, the liquid REV will be much smaller than the solid REV, and the experimental dimensions may be based upon the porosity of the medium. However, when other factors such as a medium's biological or chemical characteristics significantly impact solute behavior, it will be necessary to reconsider the definition of a continuum. For instance, most organochlorines do not ionize appreciably in typical hydrogeological settings, and as a consequence, they are not adsorbed onto the cation exchange sites of mineral soils. Most of the chemicals from this group are, however, fat soluble and are thus readily absorbed by lipids in soil organic matter. Therefore, when establishing experimental dimensions for a study involving an organochlorine, the percentage of soil organic material will be a primary criterion.

In summary, the accurate interpretation of experimental results depends upon the integrity of the relationship between the physical and mathematical problems. In our case, describing solute behavior with continuous functions implies a physical continuum. To serve its purpose, the continuum's smallest unit must consist of an averaging volume which includes enough material to yield true statistical inferences, ideally a random population of REVs. We define such a fundamental unit as a sphere with diameter $\bar{d}_{REV}$.

### 3.2. GOVERNING EQUATION

The cylindrical column shown in figure 3.1 is used to represent any porous medium with a constant cross-section greater than $\bar{d}_{REV}$. Consequently, the permeable solid will be homogeneous so that the overall cross-sectional area $(A)$ will reflect average porosity $(n)$. Flow will occur only through the void space $(An)$ and the column will have the fixed length $(\ell)$.

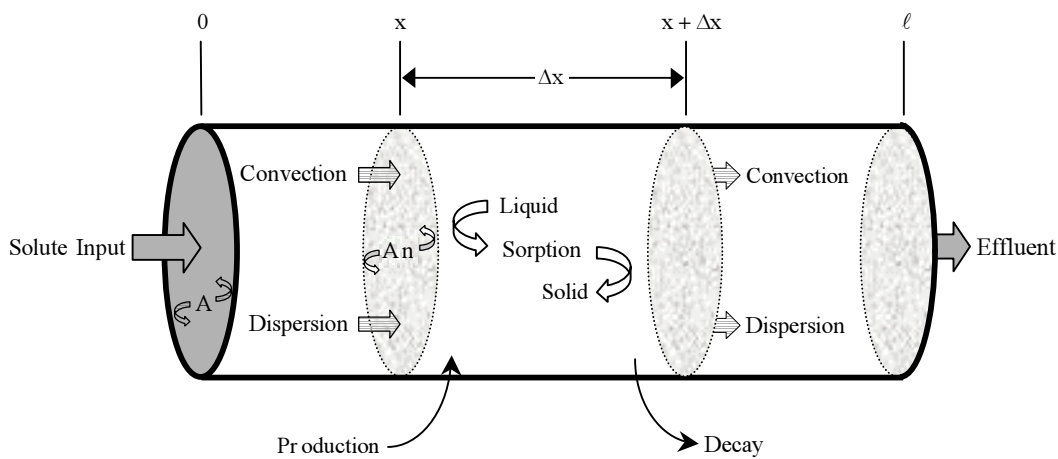

FIGURE 3.1. COLUMN SCHEMATIC.

The solid matrix is assumed to be continuously saturated with a homogeneous liquid consisting of a solvent carrying a single miscible solute whose initial concentration may vary longitudinally. A similar liquid with a concentration that may vary with time is then applied evenly over the inflow face at $x = 0$, which admits the assumption of one-dimensional transport (i.e., convection and dispersion occur only orthogonal to planes parallel with the inflow face). An element of length $\Delta x$ denotes the distance between any two such planes and defines a section with total volume $A\Delta x$ and interstitial volume $An\Delta x$.

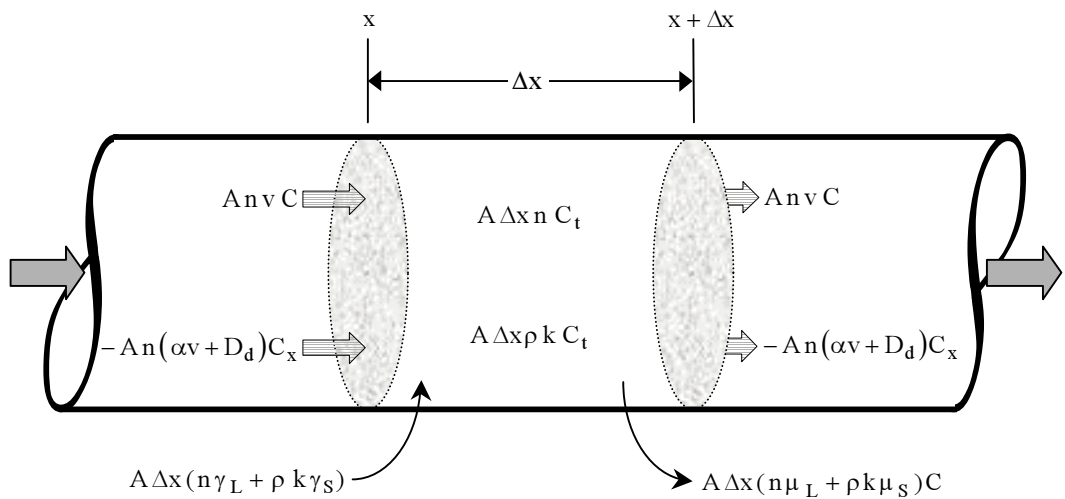

FIGURE 3.2. INTERIOR SECTION.

   The conservation of mass for the interior element shown in figure 3.2 yields the CDE that governs the physical behavior to be modeled. The rate of change of solute mass within the section can be described as

$$\text{Accumulation} = A \Delta x (n + \rho k) C_t,$$

where the constants $\rho$ and n denote the porous-medium bulk density and intrinsic porosity, respectively, and k is the linear distribution-coefficient. Transformations which result in the apparent appearance or disappearance of the solute species will be assumed to take place as

$$\text{Production} = A \Delta x \left( n \gamma_L + \rho k \gamma_S \right),$$

$$\text{Decay} = A \Delta x \left( n \mu_L + \rho k \mu_S \right) C,$$

where the constants $\gamma, \mu$ describe zero-order production and first-order decay, while the subscripts indicate whether the reaction occurs in the liquid (L) or sorbed (S) phase. Mass enters the section at x and exits at $x + \Delta x$ at the same rate of flow,

$$\text{Inflow} = \text{Outflow} = An \left( vC - (\alpha v + D_d) C_x \right).$$

Here, $\alpha$ is dispersivity, a macroscopic term which describes the spread of solute caused mostly by microscopic variations in the interstitial velocity. The constant $D_d$ describes the diffusive (Fickian) component of hydrodynamic dispersion. Substituting the above expressions into the familiar mass-balance statement

$$\text{Accumulation} = \text{Production} - \text{Decay} + \text{Inflow} - \text{Outflow}$$

yields the preliminary form of the governing equation

$$\begin{aligned} A \Delta x (n + \rho k) C_t &= A \Delta x \left( n \gamma_L + \rho k \gamma_S \right) - A \Delta x \left( n \mu_L + \rho k \mu_S \right) C \\ &\quad + An \left( vC(x,t) - (\alpha v + D_d) C_x(x,t) \right) - An \left( vC(x + \Delta x, t) - (\alpha v + D_d) C_x(x + \Delta x, t) \right). \end{aligned} \quad (3.1)$$

By rearranging (3.1) to provide difference quotients for the concentration and its spatial derivative, we obtain

$$\begin{aligned} \left( 1 + \rho k n^{-1} \right) C_t &= (\alpha v + D_d) \left( \frac{C_x(x + \Delta x, t) - C_x(x,t)}{\Delta x} \right) - v \left( \frac{C(x + \Delta x, t) - C(x,t)}{\Delta x} \right) \\ &\quad + \left( \gamma_L + \rho k n^{-1} \gamma_S \right) - \left( \mu_L + \rho k n^{-1} \mu_S \right) C, \quad n \neq 0. \end{aligned} \quad (3.2)$$

Then, for $\Delta x$ sufficiently small, (3.2) will yield

$$(1+\rho k n^{-1})C_t = (\alpha v + D_d)C_{xx} - vC_x - (\mu_L + \rho k n^{-1}\mu_S)C + (\gamma_L + \rho k n^{-1}\gamma_S). \qquad (3.3)$$

Although only the definition of a derivative is required to justify the mathematical passage from (3.2) to (3.3), the infimum of $\Delta x$ has an important meaning in the physical problem. That is, to satisfy the requirements for a material continuum, $\Delta x$ may never be smaller than $\overline{d}_{REV}$, the shortest distance over which a change in properties may be observed. The mathematical hypothesis that $\Delta x$ is sufficiently small is then equivalent to the physical requirement that $\overline{d}_{REV}$ be much less than $\ell$ which, in turn, is equivalent to requiring physical inhomogeneities in the liquid or solid medium to be compensated for by a corresponding increase in the length of the experimental flow domain. The statement that (3.3) is justified then implies that the assumption of macroscopic homogeneity has been satisfied. We may thus define the (non-negative) constants

$$\begin{aligned}\gamma &= \gamma_L + \rho k n^{-1}\gamma_S, \\ \mu &= \mu_L + \rho k n^{-1}\mu_S, \\ R &= 1 + \rho k n^{-1}, \\ D &= \alpha v + D_d,\end{aligned} \qquad (3.4)$$

which allows (3.3) to be written in the simpler form (1.1).

3.3. Boundary Conditions

A characterization of the boundaries as planes requires some subtle yet noteworthy assumptions on the magnitudes of production, decay, and the variation in the input concentration. The required restrictions arise quite naturally from the conservation of mass for a boundary layer adjacent to the entry plane of thickness $\Delta x$.

As illustrated in figure 3.3, the specific discharge-rate q and solute-input concentration g(t) determine the rate at which mass enters the vessel $Aqg(t)$. Using this term for inflow and recognizing that the expressions for accumulation, production, decay, and outflow are identical to those described for the interior section, the mass balance yields

$$An\Delta x RC_t = Aqg(t) - An(vC(\Delta x,t) - DC_x(\Delta x,t)) + An\Delta x\gamma - An\Delta x\mu C. \qquad (3.5)$$

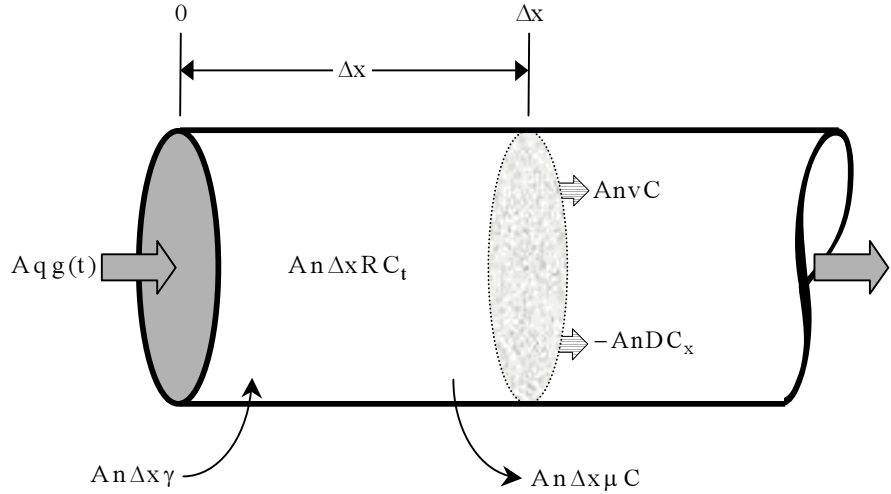

FIGURE 3.3. UPPER BOUNDARY LAYER.

Using the fact that the specific discharge-rate and interstitial velocity are equated as

$$v = qn^{-1},$$

(3.5) can be written as

$$vg(t) = \underbrace{vC(\Delta x, t) - DC_x(\Delta x, t)}_{(i)} + \underbrace{\Delta x \big( RC_t(x,t) + \mu C(x,t) - \gamma \big)}_{(ii)}, \tag{3.6}$$

from which it is clear that the boundary layer will reduce to a plane if the fate term (ii) is vanishingly small in comparison to the transport term (i). Recall that, as a condition of (3.3), $\Delta x$ must be small and, to satisfy a continuum, $\Delta x$ must have the infimum $\overline{d}_{REV}$. Since $\overline{d}_{REV}$ is defined to yield average values for $\mu, \gamma, C,$ when these are not too large, $\Delta x(\mu C(x,t) - \gamma)$ will be insignificant and (3.6) will reduce to

$$vg(t) = \underbrace{vC(\Delta x, t) - DC_x(\Delta x, t)}_{(i)} + \underbrace{\Delta x RC_t(x,t)}_{(ii)}. \tag{3.7}$$

We will obtain the desired boundary plane if we can dismiss (3.7.ii), for which we require $C_t$ to remain small at all times. That condition is satisfied by our assumption that g(t) is continuously differentiable. Thus, when we restrict our consideration to the large class of practical problems for which concentrations are dilute, transformation rates are moderate, and inputs are relatively smooth, the upper boundary condition reduces to (1.2).

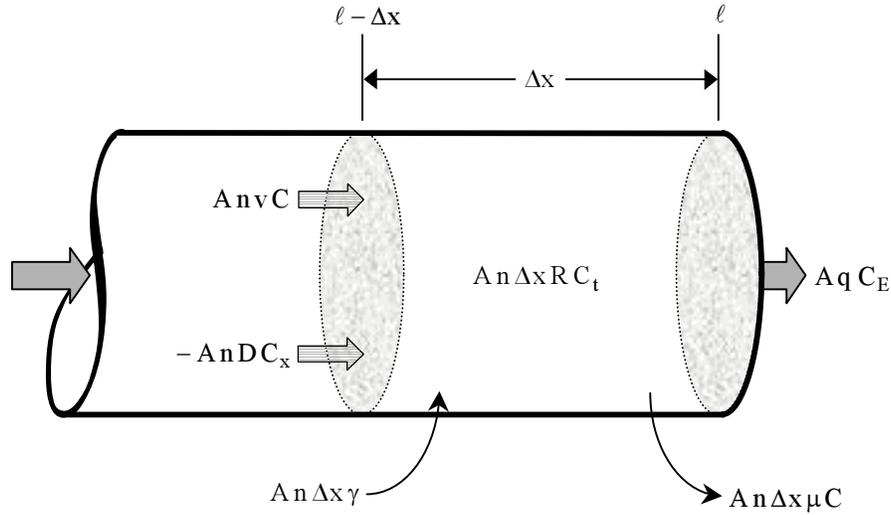

FIGURE 3.4. LOWER BOUNDARY LAYER.

Figure 3.4 illustrates mass transport through a boundary layer adjacent to the exit plane. Using the mass flow-rate $A q C_E$ for outflow, and recognizing that inflow, accumulation, production, and decay are identical to the expressions developed above, the mass balance yields

$$v C_E = \underbrace{v C(\ell - \Delta x, t) - D C_x(\ell - \Delta x, t)}_{(i)} + \underbrace{\Delta x \left( R C_t(x,t) + \mu C(x,t) - \gamma \right)}_{(ii)}. \quad (3.8)$$

If we now follow the arguments presented above, (3.8) reduces to (1.3) and thus completes the derivation of the governing equation (1.1) and the boundary conditions (1.2),(1.3).

**CHAPTER 4. PRELIMINARY TRANSFORMATIONS AND SEPARATION OF VARIABLES**

This chapter provides linear transformations which take the CDE to an inhomogeneous diffusion equation and the boundary conditions to a homogeneous form. The variables are then formally separated in the context of an outline for a justification of the separated solution.

4.1. TRANSFORMATION TO AN INHOMOGENEOUS DIFFUSION EQUATION

This section follows techniques from ordinary differential equations (ODEs) to derive a novel approach to transforming (1.1) into an inhomogeneous diffusion equation. We begin with a result referred to as reduction of order, which states that when one solution of a linear second-order ODE has been determined a second linearly independent solution can be defined from the product of two functions (see, e.g., Boyce and DiPrima 1986, sec. 3.4). In our case, it is easy to see that we may let

$$C(x,t;r,s) = u(x,t)\,\xi(x,t;r,s), \qquad (4.1)$$

if $u, \xi$ are defined in such a way that their product constitutes a linear combination; then, when (4.1) and its partial derivatives

$$C_t = u_t\,\xi + u\,\xi_t,$$

$$C_x = u_x\,\xi + u\,\xi_x,$$

$$C_{xx} = u_{xx}\,\xi + 2u_x\,\xi_x + u\,\xi_{xx}$$

are substituted in (1.1), we will have

$$R(u_t\,\xi + u\,\xi_t) = D(u_{xx}\,\xi + 2u_x\,\xi_x + u\,\xi_{xx}) - v(u_x\,\xi + u\,\xi_x) - \mu(u\,\xi) + \gamma,$$

or its equivalent

$$R\xi u_t - D\xi u_{xx} - \gamma = (2D\xi_x - v\xi)u_x - (R\xi_t - D\xi_{xx} + v\xi_x + \mu\xi)u. \qquad (4.2)$$

If we allow u to be any smooth function, which will remain free to provide a solution to the original problem, then $\xi$ may be any appropriately defined particular solution. By defining

$$\xi(x,t;r,s) = e^{rx - st}, \qquad (4.3)$$

we can follow a well-known method for solving homogeneous differential equations with constant coefficients (see, e.g., Boyce and DiPrima 1986, sec. 3.5) which will allow us to reduce the parenthetical expressions in (4.2) to constant factors. To begin, note that the partial derivatives

$$\xi_t = -s e^{rx-st},$$

$$\xi_x = r e^{rx-st},$$

$$\xi_{xx} = r^2 e^{rx-st},$$

allow us to rewrite (4.2) as

$$R\xi u_t - D\xi u_{xx} - \gamma = \xi(2Dr - v)u_x + \xi(Rs + Dr^2 - vr - \mu)u,$$

or since $\xi$ is nowhere zero,

$$R u_t - D u_{xx} - \gamma \xi^{-1} = \underbrace{(2Dr - v)}_{(i)} u_x + \underbrace{(Rs + Dr^2 - vr - \mu)}_{(ii)} u. \tag{4.4}$$

The convection term $u_x$ will now disappear when we define

$$r = \frac{v}{2D}, \quad D \neq 0,$$

which forces the coefficient (4.4.i) to zero. To eliminate u, we choose s as

$$s = \frac{1}{R}\left(\frac{v^2}{4D} + \mu\right), \quad R \neq 0,$$

so that (4.4.ii) also vanishes.

We now have the inhomogeneous diffusion equation

$$u_t = \frac{D}{R} u_{xx} + \frac{\gamma}{R} \xi^{-1}(x, t), \tag{4.5}$$

and substituting the right-hand side (rhs) of (4.1) into (1.4), we obtain the transformed initial condition

$$u(x, t_0) = \xi^{-1}(x, t_0)\phi(x), \tag{4.6}$$

and that substitution in (1.2),(1.3) yields the boundary conditions

$$u(0,t) - \frac{1}{r} u_x(0,t) = 2\xi^{-1}(0,t) g(t), \quad r \neq 0,$$

$$u(\ell,t) - \frac{1}{r} u_x(\ell,t) = 2\xi^{-1}(\ell,t) C_E. \tag{4.7}$$

## 4.2. TRANSFORMATION TO A SYSTEM WITH HOMOGENEOUS BOUNDARY CONDITIONS

We now wish to introduce a change of variable that will transform the boundary conditions on u into a homogeneous form. We thus define

$$u(x,t) = w(x,t) + e^{st} H(x,t), \tag{4.8}$$

and when this is substituted for u in the boundary conditions, we get

$$w(0,t) - \frac{1}{r} w_x(0,t) = 2\xi^{-1}(0,t) g(t) - e^{st}\left(H(0,t) - \frac{1}{r} H_x(0,t)\right),$$

$$w(\ell,t) - \frac{1}{r} w_x(\ell,t) = 2\xi^{-1}(\ell,t) C_E - e^{st}\left(H(\ell,t) - \frac{1}{r} H_x(\ell,t)\right),$$

whose rhs will be zero if H can be chosen so that $e^{st}\left(H(0,t) - (1/r) H_x(0,t)\right) = 2\xi^{-1}(0,t) g(t)$ and $e^{st}\left(H(\ell,t) - (1/r) H_x(\ell,t)\right) = 2\xi^{-1}(\ell,t) C_E$. The form of H suggests a definition in periodic functions, and it is easy to verify that

$$H(x,t) = \left(1 + \cos\frac{\pi x}{\ell}\right) g + \left(1 - \cos\frac{\pi x}{\ell}\right) e^{-r\ell} C_E \tag{4.9}$$

yields the desired result. That is, since $e^{st} H(0,t) = 2\xi^{-1}(0,t) g(t)$ and $e^{st} H(\ell,t) = \xi^{-1}(\ell,t) C_E$, we only require $H_x(0,t) = H_x(\ell,t) = 0$, but that is satisfied by

$$H_x(x,t) = \frac{\pi}{\ell}\left(e^{-r\ell} C_E - g\right) \sin\frac{\pi x}{\ell}.$$

To arrive at the differential equation in $w(x,t)$, note that substituting the rhs of (4.8) for $u(x,t)$ in (4.5) and carrying out the indicated differentiation will yield

$$w_t = \frac{D}{R} w_{xx} + e^{st}\left(\frac{\gamma}{R} e^{-rx} - (sH + H_t) + \frac{D}{R} H_{xx}\right), \tag{4.10}$$

where

$$H_t = \left(1 + \cos\frac{\pi x}{\ell}\right)\dot{g} + \left(1 - \cos\frac{\pi x}{\ell}\right)e^{-r\ell}\dot{C}_E,$$

with $\dot{g}, \dot{C}_E$ denoting t derivatives, and

$$H_{xx} = \frac{\pi^2}{\ell^2}\left(e^{-r\ell}C_E - g\right)\cos\frac{\pi x}{\ell}.$$

To facilitate an inspection of exit conditions, it is convenient to rewrite (4.10) in the equivalent form

$$w_t = \frac{D}{R}w_{xx} + e^{st}F(x,t), \qquad F(x,t) = (F_1 + F_2), \tag{4.11}$$

where

$$F_1 = \frac{\gamma}{R}e^{-rx} - \left(\left(\frac{\pi^2 D}{\ell^2 R} + s\right)\cos\frac{\pi x}{\ell} + s\right)g - \left(1 + \cos\frac{\pi x}{\ell}\right)\dot{g},$$

$$F_2 = \left(\left(\frac{\pi^2 D}{\ell^2 R} + s\right)\cos\frac{\pi x}{\ell} - s\right)e^{-r\ell}C_E - \left(1 - \cos\frac{\pi x}{\ell}\right)e^{-r\ell}\dot{C}_E.$$

Note that $F(x,t)$ simply parses the inhomogeneous term $(\gamma/R)e^{-rx} - (sH + H_t) + (D/R)H_{xx}$ so that $F_1$ is an expression in the input concentration and $F_2$ in the exit concentration.

If we now substitute the rhs of (4.8) into the auxiliary conditions for u, we will have the initial condition

$$w(x,t_0) = e^{st_0}\left(e^{-rx}\phi(x) - H(x,t_0)\right), \tag{4.12}$$

with the boundary conditions

$$w_x(0,t) - r w(0,t) = 0,$$
$$w_x(\ell,t) - r w(\ell,t) = 0. \tag{4.13}$$

4.3. Separation of Variables

In this section, we formally separate variables and provide an outline of the elements which will be required to complete our *a posteriori* justification of the separated solution

$$w(x,t) = \varphi(x)T(t). \tag{4.14}$$

Suppose that (4.11) admits the solution (4.14) so that, setting $F(x,t) = 0$, we obtain the homogeneous differential equation

$$\frac{\partial(\varphi T)}{\partial t} = \frac{D}{R}\frac{\partial^2(\varphi T)}{\partial x^2}. \tag{4.15}$$

Since we are not interested in the trivial solution $w(x,t) \equiv 0$ for all $x, t$, we may divide by $\varphi T$ to get

$$\frac{\varphi''}{\varphi} = \frac{R}{D}\frac{\dot{T}}{T}. \tag{4.16}$$

Assume that $\varphi, T$ are sufficiently smooth so that we may differentiate once more with respect to x, which yields $(\varphi''/\varphi)' = 0$. Thus, the rhs must equal a constant, call it $-\lambda$, from which we obtain the spatial problem

$$\varphi'' = -\lambda\varphi. \tag{4.17}$$

We can now show that (4.17) will have solutions which are continuously defined and differentiable any number of times:

LEMMA (ANALYTICITY): For fixed $\{\lambda: 0 < \lambda < \infty\}$,[1] the ordinary differential equation (4.17) will have Taylor (Maclaurin) series solutions which will converge absolutely for any fixed $\{x: 0 < x < \infty\}$ and which will converge uniformly to continuous functions in any closed interval $\{[x_0, x_\ell]: -\infty < x_0 \leq x \leq x_\ell < \infty\}$.

PROOF: Let $\Phi_n(x;\lambda) = \sum_n a_n x^n$ solve $\Phi'' = -\lambda\Phi$ so that $\sum_n ((n+1)(n+2)a_{n+2} + \lambda a_n) x^n = 0$. Upon dismissing $x \equiv 0$, we find the recursion formula $a_{n+2} = -\lambda a_n / ((n+1)(n+2))$. Our solutions will be defined differently for: (i) $\lambda = 0$, (ii) $\lambda > 0$, (iii) $\lambda < 0$:

(i)  For $\lambda = 0$, we must have $a_n = 0$ for all $n = 1, 2, 3, \ldots, N, \ldots$ which implies the trivial solution $\Phi \equiv 0$.

(ii) For $\lambda = \eta > 0$, let $a_n = a_n^\eta$ and for $k = 0, 1, 2, \ldots$, the recursion relation yields series that differ according to whether: (a) $n = 2k$, (b) $n = 2k+1$:

(a) The index $n = 2k$ yields a series in even powers of x,

---

[1] We will be interested in solutions for point values of λ only. However, uniform convergence for a continuous λ could easily be provided by interchanging the symbols x and λ in each of the proposition sets so that λ would have terminal values. For a discussion of uniform convergence, see Gustafson (1999, sec. 1.6.3, theorem 2).

$$\Phi_{n,0}(x;\eta) = a_0^\eta \left( 1 - \frac{(x\sqrt{\eta})^2}{2!} + \frac{(x\sqrt{\eta})^4}{4!} - \frac{(x\sqrt{\eta})^6}{6!} + \cdots + \frac{(-1)^n (x\sqrt{\eta})^{2n}}{(2n)!} + \cdots \right),$$

which converges absolutely in the whole $x, \lambda$ plane with its d'Alembert ratio $|\Phi_{n+1,0}/\Phi_{n,0}| \to 0$ as $n \to \infty$. Since $\sum_1^N \Phi_{n,0}(x;\eta)$ is a partial sum of the Maclaurin series for $a_0^\eta \text{Cos}(x\sqrt{\eta})$, we may fix $\eta$, and for any $\varepsilon$, however small, there will be an $N$ for which $|a_0^\eta \text{Cos}(x\sqrt{\eta}) - \sum_1^N \Phi_{n,0}(x;\eta)| < \varepsilon$ in any interval $[x_0, x_\ell]$ such that the convergence is uniform.

(b) The index $n = 2k + 1$ yields the series in odd powers of $x$

$$\Phi_{n,1}(x;\eta) = a_1^\eta \left( x\sqrt{\eta} - \frac{(x\sqrt{\eta})^3}{3!} + \frac{(x\sqrt{\eta})^5}{5!} - \frac{(x\sqrt{\eta})^7}{7!} + \cdots + \frac{(-1)^n (x\sqrt{\eta})^{2n+1}}{(2n+1)!} + \cdots \right),$$

which again converges absolutely. Uniform convergence then follows as above from the fact that $\sum_1^N \Phi_{n,1}(x;\eta)$ is a partial sum for $a_1^\eta \text{Sin}(x\sqrt{\eta})$.

Taking (a) and (b) together, $\sum_1^N (\Phi_{n,0}(x;\eta) + \Phi_{n,1}(x;\eta))$ will converge absolutely and uniformly to the continuous function $\Phi(x;\eta) = a_0^\eta \text{Cos}(x\sqrt{\eta}) + a_1^\eta \text{Sin}(x\sqrt{\eta})$ in $[x_0, x_\ell]$.

(iii) For $\lambda = \alpha < 0$, letting $a_n = a_n^\alpha$ we again have series which differ according to whether: (a) $n = 2k$, (b) $n = 2k + 1$:

(a) For $n = 2k$,

$$\Phi_{n,0}(x;\alpha) = a_0^\alpha \left( 1 + \frac{(x\sqrt{\alpha})^2}{2!} + \frac{(x\sqrt{\alpha})^4}{4!} + \frac{(x\sqrt{\alpha})^6}{6!} + \cdots + \frac{(x\sqrt{\alpha})^{2n}}{(2n)!} + \cdots \right),$$

which converges absolutely for all $x, \lambda$ and uniformly to the continuous hyperbolic function $a_0^\alpha \text{Cosh}(x\sqrt{\alpha})$ in $[x_0, x_\ell]$.

(b) For $n = 2k + 1$,

$$\Phi_{n,1}(x;\alpha) = a_1^\alpha \left( x\sqrt{\alpha} - \frac{(x\sqrt{\alpha})^3}{3!} + \frac{(x\sqrt{\alpha})^5}{5!} - \frac{(x\sqrt{\alpha})^7}{7!} + \cdots + \frac{(x\sqrt{\alpha})^{2n+1}}{(2n+1)!} + \cdots \right),$$

which converges absolutely for all $x, \lambda$ and uniformly to $a_1^\alpha \text{Sinh}(x\sqrt{\alpha})$ in $[x_0, x_\ell]$. Taking (a) and (b) together, $\sum_1^N (\Phi_{n,0}(x;\alpha) + \Phi_{n,1}(x;\alpha))$ will converge absolutely and uniformly to the continuous function $\Phi(x;\alpha) = a_0^\alpha \text{Cosh}(x\sqrt{\alpha}) + a_1^\alpha \text{Sinh}(x\sqrt{\alpha})$ in $[x_0, x_\ell]$.

We now resume our separation with the substitution of (4.14) into the boundary conditions (4.13), which yields

$$T(t)\big(\varphi'(0) - r\varphi(0)\big) = 0,$$

$$T(t)\big(\varphi'(\ell) - r\varphi(\ell)\big) = 0, \quad (4.18)$$

which will be satisfied for any bounded $T(t)$ (i.e., there will be a $K > 0$ such that $|T(t)| \leq K$ for all $\{t : t_0 \leq t\}$).[2] Then,

$$\varphi'(0) - r\varphi(0) = 0,$$

$$\varphi'(\ell) - r\varphi(\ell) = 0 \quad (4.19)$$

will be satisfied for some choice of the constants $a_0^\eta, a_1^\eta, \eta$ and $a_0^\alpha, a_1^\alpha, \alpha$.

If we now suppose that

$$F(x,t) = \varphi(x)f(t) \quad (4.20)$$

furnishes a proper definition of F, we may write (4.11) as

$$\frac{\partial(\varphi T)}{\partial t} = \frac{D}{R}\frac{\partial^2(\varphi T)}{\partial x^2} + e^{st}\varphi f \quad (4.21)$$

which, upon division by $\varphi T$, yields

$$\frac{\varphi''}{\varphi} = \frac{R}{D}\left(\frac{\dot{T}}{T} - \frac{f}{T}\right). \quad (4.22)$$

Then, differentiating (4.22) once more with respect to x will yield $(\varphi''/\varphi)' = 0$, implying that $\varphi''/\varphi$ is again equal to some constant. Therefore, the spatial solutions $\varphi$ of the inhomogeneous problem (4.21),(4.22) can be furnished from a solution of the homogenous problem (4.15),(4.16) provided that $F(x,t)$ admits an expansion like (4.20).

The justification for our separated solution will be provided over the next three chapters. Chapter 5 will show that the solutions of the spatial system (4.17),(4.19) compose an infinite linear combination whose finite subsets have linear bases. In chapter 6, we will furnish a solution for the differential equation (4.21) that satisfies (4.12),(4.18), and we will go on to show that our solution converges uniformly and unconditionally to its defining function, thus verifying the expansion (4.20). Chapter 7 then uses the demonstrated convergence to deduce that (4.11)-(4.13) has a uniform approximation in a linear basis.

---

[2]As part of the problem statement in section 1.2, we asked for a time interval in which the input concentration would be continuously differentiable. Then, when discussing the upper boundary in sec. 3.3, we explicitly excluded highly discontinuous inputs and large transformation rates. It is thus reasonable to suppose that $T(t)$ will be continuously differentiable for $(t_0, t)$ and bounded for $[t_0, t]$.

**CHAPTER 5. STURM-LIOUVILLE ROBIN**

5.1. DEFINITIONS AND PROPERTIES OF THE OPERATORS

DEFINITION 5.1 (STURM-LIOUVILLE OPERATOR): We define the Sturm-Liouville operator[1]

$$S(\varphi(x)) := \varphi'' = -\lambda \varphi, \quad 0 < x < \ell$$

from the spatial differential equation (4.17).

We may assume that any solution $\varphi$ of $S(\varphi)$ is real valued and that $\lambda$ is a real number. This will not involve any loss of generality since we show in the sequel that any solution of S must be entirely real. Solutions of S then attain the familiar properties of a real-valued linear operator.[2]

LEMMA 5.2 (LINEARITY): The operator $S(y)$ that takes y to the function $y'' + \lambda y$ identically zero in the open interval $(0, \ell)$ constitutes a linear transformation.

PROOF: From the definition of linearity, S will be linear if, for $c_1, c_2$ any two numbers and $y_1, y_2$ any two solutions of $S(y)$, $S(c_1 y_1 + c_2 y_2) = c_1 S(y_1) + c_2 S(y_2)$, which follows immediately, since

$$S(c_1 y_1 + c_2 y_2) = (c_1 y_1 + c_2 y_2)'' + \lambda (c_1 y_1 + c_2 y_2)$$

$$= c_1 y_1'' + c_2 y_2'' + \lambda c_1 y_1 + \lambda c_2 y_2$$

$$= c_1 (y_1'' + \lambda y_1) + c_2 (y_2'' + \lambda y_2)$$

$$= c_1 S(y_1) + c_2 S(y_2)$$

$$= 0.$$

DEFINITION 5.3 (ROBIN BOUNDARY OPERATORS): The Robin boundary operators

$$B_1(\varphi) := \varphi'(0) - r\varphi(0) = 0,$$

$$B_2(\varphi) := \varphi'(\ell) - r\varphi(\ell) = 0,$$

are defined by the conditions (4.19).

---

[1] For alternative discussions, cf. (i) Sturm-Liouville problems in the context of the solution of partial differential equations in Strauss (1992, sec. 11.4) or (ii) the much broader Sturmian theory of ordinary differential equations in Ince (1956, chap. 10).

[2] For alternative discussions, cf. (i) the properties of linear operators in Berg and McGregor (1966, sec. 1.3) or (ii) linear differential operators in Boyce and DiPrima (1986, theorem 3.2).

DEFINITION 5.4 (SLR): The Sturm-Liouville operator S together with the Robin boundary conditions $B_1, B_2$ composes a system which we will refer to concisely as Sturm-Liouville Robin (SLR).

A function y that satisfies $S(y)$ alone will be referred to simply as a solution of S. Such solutions will satisfy the boundary conditions only for carefully selected numbers $\beta$. Since there are infinitely many such $\beta$, they can be regarded as the characteristic numbers which, in turn, establish characteristic functions $y(x;\beta)$, which are also infinite in number. We define these characteristic solutions as follows:

DEFINITION 5.5 (EIGENSOLUTIONS): We refer to a solution y of $S(y)$ as an eigenfunction and a number $\lambda$ as an eigenvalue only when $y(x;\lambda)$ satisfies $B_1(y)$ and $B_2(y)$ (i.e., only when y completely solves SLR).

We will provide the eigensolutions of SLR by establishing some general properties and then steadily increasing the level of specification until the eigensolutions emerge as particular consequences. To begin, consider LaGrange's identity

$$y_2 S(y_1) - y_1 S(y_2) = y_1'' y_2 - y_1 y_2'' \tag{5.1}$$

and its secondary identity

$$y_1'' y_2 - y_1 y_2'' = (y_1' y_2 - y_1 y_2')'. \tag{5.2}$$

Note that if we restrict definition 5.1 to $S(y) := y''$, then LaGrange's identities are valid for any two twice continuously differentiable functions. If, however, we require $S(y)$ to be identically satisfied (i.e., $y_1'' = -\beta_1 y_1$, $y_2'' = -\beta_2 y_2$), LaGrange's identity will yield

$$(\beta_2 - \beta_1) y_1 y_2 = y_1'' y_2 - y_1 y_2''. \tag{5.3}$$

DEFINITION 5.6 (WRONSKIAN): Let $y_1(x), y_2(x)$ be any two solutions of the same linear second-order differential equation. We then have

$$W(y_1, y_2)(x) = y_1' y_2 - y_1' y_2,$$

the Wronskian of $y_1, y_2$.

LEMMA 5.7 (WRONSKIAN): Let $y_1(x;\beta), y_2(x;\beta)$ be any two eigensolutions of $S(y)$. Then, $W(y_1, y_2)(x)$ will be constant in $0 < x < \ell$.

PROOF: Equating (5.2) to (5.3) yields $(y_1' y_2 - y_1 y_2')' = 0$ whenever $\beta_2 = \beta_1$, and since $y_1, y_2$ are solutions of $S(y)$, this holds throughout $(0, \ell)$. Then, from the definition of the Wronskian, $(y_1' y_2 - y_1 y_2')' = W'(y_1, y_2)(x)$.

Thus far, it has sufficed to ask that the solutions of S be defined in the open interval $(0, \ell)$. As we continue, we will require the integral of various equalities between (5.1)-(5.3). This involves the assumption that the integrands are bounded on the closed interval $[0, \ell]$, so we provide for that with the following comment:

COMMENT 5.8 (INTEGRABILITY): Let $\mathcal{R}(x)$, which may be a product of other functions, be continuously defined on the open interval $x_0 < x < x_\ell$ with no more than a finite discontinuity at its end points $x_0, x_\ell$. Then we say that $\mathcal{R}(x)$ is Riemann integrable for the open interval $(x_0, x_\ell)$.

Equating (5.2) to (5.3) and integrating over $(0, \ell)$, we obtain Green's (second) formula

$$(\beta_2 - \beta_1) \int_0^\ell y_1 y_2 \, dx = \int_0^\ell (y_1'' y_2 - y_1 y_2'') \, dx = y_1' y_2 - y_1 y_2' \Big|_0^\ell, \quad (5.4)$$

which, as a consequence of comment 5.8, will be valid for $y_1(x; \beta_1), y_2(x; \beta_2)$ any two solutions of $S(y)$ irrespective of boundary conditions.

When we require each of $y_1, y_2$ to satisfy both of the boundary conditions $B_1, B_2$, they become eigenfunctions and will therefore fulfill the symmetry conditions

$$y_1 B_1(y_2) - y_2 B_1(y_1) = y_1' y_2 - y_1 y_2' \Big|_{x=0} = 0,$$

$$y_2 B_2(y_1) - y_1 B_2(y_2) = y_1' y_2 - y_1 y_2' \Big|_{x=\ell} = 0. \quad (5.5)$$

Note that (5.5) can also be obtained from a Cramer's rule solution of the algebraic systems $B_i(y_1), B_i(y_2)$ where $i = 1, 2$. That is, when $y_1, y_2$ are taken as coefficients in those systems, if there is to be a non-trivial solution, the coefficient matrix must vanish, which again yields the rightmost equalities of (5.5).

LEMMA 5.9 (EIGENVALUE MULTIPLICITY): SLR may have no more than one unique eigenfunction for any fixed eigenvalue.[3]

---

[3] Alternate proofs of lemma 5.9 can be provided by noting that, because of the symmetry, the vanishing of W can be deduced from: (i) Rolle's theorem (see, e.g., Widder 1989, sec. 2.3); (ii) the law of the mean (see, e.g., Widder 1989, sec. 2.4); or (iii) a slight reformulation of lemma 5.17 (see footnote 5).

PROOF: Let $y_1, y_2$ be eigenfunctions belonging to the same eigenvalue $\beta$. Then, $y_1'' y_2 - y_1 y_2''$ $= (y_1' y_2 - y_1 y_2')'$ will be bounded on $[0, \ell]$ and the integral

$$\int_0^\ell (y_1' y_2 - y_1 y_2')' \, dx = (y_1' y_2 - y_1 y_2') \Big|_0^\ell$$

is well defined. Since $y_1, y_2$ each satisfy $B_1, B_2$, symmetry guarantees the vanishing of the rhs, which is precisely $W(y_1, y_2)(x)$. Then, since $W$ is independent of x, it must vanish identically in $(0, \ell)$ so that $y_1, y_2$ are linearly dependent.[4]

LEMMA 5.10 (REAL EIGENSOLUTIONS):

1. Any eigenvalue will be a real number;
2. The associated eigenfunctions may be chosen to be real valued.

PROOF:

1. For $i = \sqrt{-1}$, let $\lambda = \mu + i\nu$ be a complex eigenvalue which defines the complex eigenfunction $\psi(x; \lambda)$. The complex conjugate of $\lambda$ will then be the eigenvalue $\overline{\lambda} = \mu - i\nu$ which will belong to the eigenfunction $\overline{\psi}(x; \overline{\lambda})$. Taking $y_1 = \psi$, $y_2 = \overline{\psi}$ in Green's formula, we get

$$-2i\nu \int_0^\ell |\psi|^2 \, dx = \psi' \overline{\psi} - \overline{\psi}' \psi \Big|_0^\ell,$$

where the rhs is identically zero because $\psi, \overline{\psi}$ satisfy symmetry by (5.5). Upon dismissing the trivial solution $\psi(x) \equiv 0$, we must have $\nu \equiv 0$ or, equivalently, $\lambda = \overline{\lambda} = \mu$.

2. Since $\psi(x; \mu), \overline{\psi}(x; \mu)$ belong to the same real eigenvalue, they will differ by no more than a constant factor. Then, since $\psi = \chi + i\zeta$ solves $S(\psi) = (\chi + i\zeta)'' = -\mu(\chi + i\zeta)$, which will separate into the real $\chi'' = -\mu\chi$ and imaginary $\zeta'' = -\mu\zeta$, the real eigenfunction $\chi(x; \mu)$ will suffice for a unique solution.

---

[4] For the fact that a vanishing (non-vanishing) Wronskian is a necessary and sufficient for the linear dependence (independence) of two solutions of the same ODE, see Tenenbaum and Pollard (1985, comment 64.15) or Boyce and DiPrima (1986, theorem 3.9).

5.2. EIGENSOLUTIONS

We are now in a position to provide explicit formulas for the solutions of SLR, and these solutions will differ according to the sign on the eigenvalues. We will define the eigenfunctions for the non-positive eigenvalues first and then proceed to the positive eigenvalues which will yield our infinite solution set.

For $\lambda = 0$, the general solution of $S(\varphi)$ is the straight line $\varphi(x) = \mathcal{A}x + \mathcal{B}$ which must comply with the end conditions. The upper boundary $B_1$ requires $\mathcal{B} = \mathcal{A}/r$ and accordingly the lower boundary $B_2$ yields $r\ell = 0$. Note that $r \neq 0$ as a condition of (4.7) so that the lower boundary will be satisfied only for $\ell \equiv 0$. Thus, $\lambda = 0$ is not in the solution set of SLR.

For any non-zero $\lambda$, the general solution for $S(\varphi)$ is

$$\varphi(x) = C_1 e^{x\sqrt{-\lambda}} + C_2 e^{-x\sqrt{-\lambda}}, \qquad \lambda \neq 0. \tag{5.6}$$

For $\lambda < 0$, we define the constants $\mathcal{D}_1 = (C_1 + C_2)/2$, $\mathcal{D}_2 = (C_1 - C_2)/2$ and set $\lambda = -\alpha^2 \, (\alpha > 0)$ so that (5.6) may be written as $\varphi(x) = \mathcal{D}_1 \operatorname{Cosh}(\alpha x) + \mathcal{D}_2 \operatorname{Sinh}(\alpha x)$. It is then easy to see that $B_1$ will require $\mathcal{D}_2 = (r/\alpha)\mathcal{D}_1$ and that $B_2$ will yield $\alpha = r$. Upon recalling the definitions for the hyperbolic functions, which in this case yield the expressions $\operatorname{Cosh}(rx) = (e^{rx} + e^{-rx})/2$ and $\operatorname{Sinh}(rx) = (e^{rx} - e^{-rx})/2$, it becomes clear that there is only one (real-valued) solution belonging to a single (real) number, for which we have the following:

DEFINITION 5.11 (NEGATIVE EIGENSOLUTION): We will refer to

$$\varphi_0(x) = e^{rx}, \qquad \lambda_0 = -r^2,$$

as the negative eigensolution, and references to the negative eigenfunction or the negative eigenvalue will be made clear.

We note in passing that this will not cause the solution to become arbitrarily large with increasing t, as often occurs with negative eigenvalues. It will be shown that, in particular, the final representation remains bounded for all t.

For $\lambda > 0$, let $\lambda = \kappa^2 \, (\kappa > 0)$ and (5.6) will yield $\varphi(x) = \mathcal{E}_1 e^{i\kappa x} + \mathcal{E}_2 e^{-i\kappa x}$. Then, introducing the constants $\mathcal{J}_1 = \mathcal{E}_1 + \mathcal{E}_2$, $\mathcal{J}_2 = i(\mathcal{E}_1 - \mathcal{E}_2)$, Euler's formula will provide for the equivalent expression $\varphi(x) = \mathcal{J}_1 \operatorname{Cos}(\kappa x) + \mathcal{J}_2 \operatorname{Sin}(\kappa x)$. By going on to define $\mathcal{P}_1 = \mathcal{J}_1$, $\mathcal{P}_2 = (\mathcal{J}_2)^2 = -(\mathcal{E}_1 - \mathcal{E}_2)^2$, we can then argue explicitly that $\mathcal{J}_2 \operatorname{Sin}(\kappa x)$, which has the imaginary coefficient, is interchangeable with $\mathcal{P}_2 \operatorname{Sin}(\kappa x)$, which has the real coefficient. The real-valued expression $\varphi(x) = \mathcal{P}_1 \operatorname{Cos}(\kappa x) + \mathcal{P}_2 \operatorname{Sin}(\kappa x)$ will then furnish a fundamental solution set for $S(\varphi)$, as verified by the following lemma:

LEMMA 5.12 (FUNDAMENTAL SOLUTION SET): For any fixed $\lambda = \kappa^2 \, (\kappa > 0)$, the linear combination $\mathcal{P}_1 y_1 + \mathcal{P}_2 y_2$ with $y_1 = \text{Cos}(\kappa x)$, $y_2 = \text{Sin}(\kappa x)$ will furnish a unique representation for every possible solution of $S(y)$. This is equivalent to the following two-part statement:
1. $y_1, y_2$ are linearly independent solutions of $S(y)$ on $0 < x < \ell$ and
2. $y_1, y_2$ are the maximal number of solutions of $S(y)$ for strictly positive $\kappa$.

PROOF:
1. By definition 5.6, the Wronskian of $y_1, y_2$ is $W(y_1, y_2)(x) = -\kappa \left( \text{Cos}^2(\kappa x) + \text{Sin}^2(\kappa x) \right)$, and since W is nowhere zero in the interval $0 < x < \ell$, $y_1, y_2$ are linearly independent (see footnote 4).
2. Since $S(y)$ is second order, it can have no more than two linearly independent solutions.

We resume our discussion with a consideration of the boundary conditions. For $B_1$, take $\mathcal{P}_2 = (r/\kappa)\mathcal{P}_1$, and since $\mathcal{P}_1$ is arbitrary, we may set $\mathcal{P}_1 = 1$. Then,

$$\varphi(x) = \text{Cos}(x\sqrt{\lambda}) + \frac{r}{\sqrt{\lambda}} \text{Sin}(x\sqrt{\lambda}) \tag{5.7}$$

will be a fundamental solution of $S(\varphi)$ which satisfies $B_1$ continuously for $\lambda > 0$. The positive eigenvalues can now be determined by asking that (5.7) satisfy $B_2$. Note that since $\lambda, r \neq 0$, $B_2$ will yield $(\kappa^2 + r^2)\text{Sin}(\kappa\ell) = 0$ so that $\kappa = n\pi/\ell$, which leads to:

DEFINITION 5.13 (POSITIVE EIGENSOLUTIONS): We refer to the positive eigensolutions as the eigenvalues

$$\lambda_n = \frac{n^2 \pi^2}{\ell^2}, \quad n = 1, 2, 3, \ldots, N, \ldots, \tag{5.8}$$

which can be arranged in the increasing sequence

$$\lambda_1 < \lambda_2 < \lambda_3 < \cdots < \lambda_N < \cdots, \quad \lim_{n \to \infty} \lambda_n = \infty,$$

together with their correspondingly well-ordered set of eigenfunctions

$$\varphi_n(x) = \text{Cos}(x\sqrt{\lambda_n}) + \frac{r}{\sqrt{\lambda_n}} \text{Sin}(x\sqrt{\lambda_n}). \tag{5.9}$$

Note that since $\varphi_n$ is a fundamental solution set for $S(\varphi_n)$, and satisfies the boundary conditions for each $n = 1, 2, 3, \ldots$, $\{\varphi_n\}$ are the unique solutions of SLR.

We now define the inner product, which furnishes a compact way of expressing the integral of the product of two functions:

DEFINITION 5.14 (INNER PRODUCT): Consider two functions $f(x), h(x)$ continuous on the interval $x_0 < x < x_\ell$ and bounded on $x_0 \leq x \leq x_\ell$. We then define

$$(f,h) = \int_{x_0}^{x_\ell} f(x)h(x)\,dx,$$

which we refer to as the inner product of $f, h$.

DEFINITION 5.15 (NORMALIZING CONSTANTS): Let $\varphi_n$ be an eigensolution of SLR for any fixed n; we then refer to the inner product

$$(\varphi_n, \varphi_n) = \int_0^\ell \varphi_n^2 \, dx, \quad n = 0, 1, 2, \ldots,$$

as the normalizing constant for $\varphi_n$.

We can show that the normalizing constants will be positive (real) numbers by providing a formula for computing their value.

LEMMA 5.16 (NORMALIZING CONSTANT, $\lambda < 0$): The negative eigensolution will have the normalizing constant

$$(\varphi_0, \varphi_0) = \frac{e^{2r\ell} - 1}{2r},$$

which will be a positive number.

PROOF: This follows from the straightforward integration of $\varphi_0^2$ over $(0, \ell)$, and since $r, \ell > 0$, $(\varphi_0, \varphi_0)$ will be positive number.

For the remaining $\varphi_n$, we define the normalizing constants by following a method which furnishes the ancillary illustration that our eigensolutions are continuous functions of $\lambda$ for the upper boundary:[5]

LEMMA 5.17 (NORMALIZING CONSTANTS, $\lambda > 0$): The positive eigensolutions will have the normalizing constants

$$(\varphi_n, \varphi_n) = \frac{(r^2 + \lambda_n)\ell}{2\lambda_n}, \quad n = 1, 2, 3, \ldots,$$

which will be positive numbers.

---

[5] With a slight reformulation, lemma 5.17 furnishes an alternative proof to our lemma 5.9 (see Ince 1956, sec. 10.72).

PROOF: Consider the differential equation

$$\frac{\partial \varphi''}{\partial \lambda} = -\lambda \frac{\partial \varphi}{\partial \lambda} - \varphi \tag{5.10}$$

and the boundary condition

$$B_1\left(\frac{\partial \varphi}{\partial \lambda}\right) = \frac{\partial \varphi'(0)}{\partial \lambda} - r\frac{\partial \varphi(0)}{\partial \lambda} = 0$$

which are obtained by differentiating $S(\varphi)$ and $B_1(\varphi)$ with respect to $\lambda$. If we now take $y_1 = \varphi$, $y_2 = \partial \varphi / \partial \lambda$ in Green's formula, we obtain

$$\int_0^\ell \left(\frac{\partial \varphi}{\partial \lambda} \varphi'' - \varphi \frac{\partial \varphi''}{\partial \lambda}\right) dx = \varphi' \frac{\partial \varphi}{\partial \lambda} - \varphi \frac{\partial \varphi'}{\partial \lambda} \bigg|_0^\ell . \tag{5.11}$$

Since by (5.10), $\partial \varphi / \partial \lambda = -(1/\lambda)(\partial \varphi'' / \partial \lambda + \varphi)$, and by definition 5.1, $\varphi'' = -\lambda \varphi$, the integrand on the left-hand side (lhs) of (5.11) will yield $\varphi^2$ so that

$$\int_0^\ell \varphi^2 \, dx = \varphi' \frac{\partial \varphi}{\partial \lambda} - \varphi \frac{\partial \varphi'}{\partial \lambda} \bigg|_0^\ell . \tag{5.12}$$

Now, note that (5.7) is analytic in $\lambda$ and satisfies $B_1(\varphi)$ continuously for all $0 < \lambda < \infty$ so that $B_1(\partial \varphi / \partial \lambda)$ is easily satisfied. The symmetry condition at $x = 0$ is then well defined for $y_1 = \varphi, y_2 = \partial \varphi / \partial \lambda$ and, consequently, (5.12) reduces to

$$\int_0^\ell \varphi^2 \, dx = \varphi'(\ell, \lambda) \frac{\partial \varphi(\ell, \lambda)}{\partial \lambda} - \varphi(\ell, \lambda) \frac{\partial \varphi'(\ell, \lambda)}{\partial \lambda}, \tag{5.13}$$

so our eigensolutions satisfy the upper boundary as a continuously differentiable function of $\lambda$.

To obtain the normalizing constants for SLR, we set $\varphi = \varphi_n$ a positive eigenfunction. Then, taking the indicated derivatives of $\varphi_n(x)$, evaluating the resulting expressions for $x = \ell$, and carrying out a term-for-term substitution in (5.13), we obtain

$$\int_0^\ell \varphi_n^2 \, dx$$

$$= \left(r\cos(\ell\sqrt{\lambda_n}) - \sqrt{\lambda_n}\sin(\ell\sqrt{\lambda_n})\right)\left(\frac{r\ell}{2\lambda_n}\cos(\ell\sqrt{\lambda_n}) - \frac{1}{2\sqrt{\lambda_n}}\left(\ell + \frac{r}{\lambda_n}\right)\sin(\ell\sqrt{\lambda_n})\right)$$

$$-\left(\cos(\ell\sqrt{\lambda_n}) + \frac{r}{\sqrt{\lambda_n}}\sin(\ell\sqrt{\lambda_n})\right) - \left(\frac{\ell}{2}\cos(\ell\sqrt{\lambda_n}) - \frac{1}{2\sqrt{\lambda_n}}(r\ell+1)\sin(\ell\sqrt{\lambda_n})\right).$$

Since by (5.8) $\lambda_n = n^2\pi^2/\ell^2$, the trigonometric arguments $\ell\sqrt{\lambda_n}$ reduce to $n\pi$. As a result, the terms with sine vanish leaving only the terms which involve the square of the cosine, thus yielding the normalizing constants. This completes the proof.

We now define an essential property of the Sturm-Liouville system, orthogonality, which will allow us to establish that our infinite set of real-valued eigenfunctions has linearly independent subsets.

DEFINITION 5.18 (ORTHOGONALITY): Let $f(x), h(x)$ be two functions as in definition 5.14, and we will refer to them as orthogonal whenever $(f, h) = 0$.

LEMMA 5.19 (ORTHOGONALITY): Any two eigenfunctions which belong to different eigenvalues will be orthogonal.

PROOF: Let $\varphi_m, \varphi_k$ be two eigenfunctions belonging to separate eigenvalues $\lambda_m, \lambda_k$. Then, by Green's formula,

$$(\lambda_m - \lambda_k)\int_0^\ell \varphi_m \varphi_k \, dx = \varphi_k' \varphi_m - \varphi_k \varphi_m' \Big|_0^\ell,$$

where the rhs is identically zero because any two eigensolutions will satisfy symmetry as defined by (5.5). Then, since $\lambda_m \neq \lambda_k$,

$$(\varphi_m, \varphi_k) = 0, \quad m \neq k.$$

As demonstrated by lemma 5.10, the eigenfunctions may be chosen to be real valued, and recall that we actually derived formulas for complex and real-valued solutions which were equivalent to within a constant factor. However, we will now see that, as a corollary of lemma 5.19, we must choose the real-valued eigenfunctions:

COROLLARY 5.20 (COMPLEX EIGENFUNCTIONS): Complex eigenfunctions do not preserve orthogonality.

PROOF: Let $\varphi = U + iV$ be a complex eigensolution, and since $S(\varphi)$ is real valued and linear, $\overline{\varphi} = U - iV$ must also be a solution. Then, to satisfy the orthogonality established in lemma 5.19, we must have $(\varphi, \overline{\varphi}) = 0$, but

$$\int_0^\ell \varphi \overline{\varphi}\, dx = \int_0^\ell (U^2 + V^2)\, dx,$$

which can be zero only if $U \equiv V \equiv 0$ implying $\varphi \equiv 0$, which is only trivially true because the zero function is orthogonal to every function including itself. This shows that complex eigenfunctions are not, in general, orthogonal.[6]

Up to this point, we have only been interested in the properties of single or paired eigensolutions. We now wish to infer that certain linear relationships exist among all members of the infinite set of solutions of SLR. Since linearity is an algebraic notion and is therefore defined only for sequences with a finite number of members, we provide for representative subsets with the following comment:

COMMENT 5.21 (COUNTABLE SUBSETS): The eigensolutions form a countably infinite set. That is, although the sequence of eigenvalues $\{\lambda_n\}_0^N$ and the associated eigenfunctions $\{\varphi_n\}_0^N$ are defined for $N \to \infty$, the index n for each member of the sequence corresponds to a single positive integer, so it will make sense to discuss algebraic properties for some arbitrarily large subset of the eigensolutions.

Lemma 5.2 will now have the following corollary:

COROLLARY 5.22 (LINEAR COMBINATIONS): A linear combination of any finite subset of the eigenfunctions is also a solution for S.

PROOF: For any constants $\{c_n\}_0^N$, define the linear combination $\widetilde{\varphi} = c_0\varphi_0 + c_1\varphi_1 + \cdots + c_N\varphi_N$ so that $S(\widetilde{\varphi}) = (c_0\varphi_0 + c_1\varphi_1 + \cdots + c_N\varphi_N)'' + \lambda(c_0\varphi_0 + c_1\varphi_1 + \cdots + c_N\varphi_N)$. Then, by induction on lemma 5.2, $\widetilde{\varphi}$ is clearly a solution of S since $S(\widetilde{\varphi}) = c_0 S(\varphi_0) + c_1 S(\varphi_1) + \cdots + c_N S(\varphi_N) = 0$.

THEOREM 5.23 (LINEAR INDEPENDENCE): A linear combination of any subset of the eigenfunctions will compose a linearly independent set for S.

---

[6] For an alternative discussion of the impossibility of complex eigenfunctions for problems like SLR, cf. Ince (1956, sec. 10.7) whose proof entails the fact that the eigenvalues are also real, providing an alternative to lemma 5.10.

PROOF: We must show that the only solution of the equation $c_0\varphi_0 + c_1\varphi_1 + \cdots + c_N\varphi_N = 0$ is $c_n = 0$ for all $n = 0, 1, \ldots, N$. To see that this holds, consider the inner product $(\varphi_m, c_0\varphi_0 + c_1\varphi_1 + \cdots + c_m\varphi_m + \cdots + c_N\varphi_N) = c_n(\varphi_m, \varphi_0 + \varphi_1 + \cdots + \varphi_m + \cdots + \varphi_N)$ where, by orthogonality, $(\varphi_m, \varphi_n) = 0$ for all $n \neq m$ so that only $c_n(\varphi_m, \varphi_m)$ remains. However, $(\varphi_m, \varphi_m) > 0$, so we must have $c_n \equiv 0$ for all n no matter how we choose N.

This chapter reduces to the essential residue of a finite sequence of eigenfunctions. Recall that there is only one negative eigensolution and that the positive eigensolutions compose a fundamental set for each $n = 1, 2, 3, \ldots$. Then, since any finite subset of the eigenfunctions will be linearly independent, $\{\varphi_n\}_0^N$ will include every unique solution of SLR up to the $N^{th}$ eigenfunction. Any N-dimensional subset of the eigenfunctions will therefore have a Hamel basis.[7]

---

[7] A set of solutions $\{y_n\}_0^N$ of a linear operator L is said to compose an N-dimensional Hamel basis for L if: (i) the $\{y_n\}_0^N$ form a linearly independent set whose (ii) N members represent every possible solution of L. Note that the fundamental solution set of lemma 5.12 is an example of a Hamel basis where $N = 2$.

**CHAPTER 6. REPRESENTATION OF THE SOLUTION**

In this chapter, we derive a representation for a solution to the original problem in C from the arguments required to furnish a justification of the separation of variables.

6.1. SOLUTION TO THE SYSTEM WITH HOMOGENEOUS BOUNDARY CONDITIONS

We will show that there is a well-behaved solution to the separated system in w by demonstrating that the terms

$$w_n(x,t) = T_n(t)\,\varphi_n(x), \quad n = 0,1,2,\ldots \tag{6.1}$$

compose a convergent series that satisfies the boundary and initial conditions continuously from the interior of the $x,t$ domain.[1] That continuity condition for the separated boundary conditions will be

$$T_n(t)\Big(\lim_{\substack{x\to 0\\x>0}} \varphi'_n(x) - r\lim_{\substack{x\to 0\\x>0}} \varphi_n(x)\Big) = 0,$$

$$T_n(t)\Big(\lim_{\substack{x\to \ell\\x<\ell}} \varphi'_n(x) - r\lim_{\substack{x\to \ell\\x<\ell}} \varphi_n(x)\Big) = 0 \tag{6.2}$$

for each n. Since each $\varphi_n$ is continuously differentiable any number of times, it should be easy to see that the limits in (6.2) are satisfied (i.e., the trigonometric functions $\varphi_n, \varphi'_n$ are everywhere continuous and must therefore approach their boundary values continuously).

It then remains for us to define the $T_n(t)$ and show that they are bounded. Substituting (6.1) into (4.11) will yield

$$\dot{T}_n\,\varphi_n = \frac{D}{R} T_n\,\varphi''_n + e^{st}\varphi_n f_n(t), \tag{6.3}$$

upon recalling our assumption that $F(x,t)$ will have an expansion like (4.20). Since the eigenpairs satisfy $\varphi''_n = -\lambda_n \varphi_n$ for each n, (6.3) may be written as $(\dot{T}_n + (D/R)\lambda_n T_n - e^{st} f_n)\,\varphi_n = 0$. Then, because $\varphi_n$ does not vanish identically for even one n, the term in parentheses must be zero for each n, and this yields the first-order ODE

---

[1] For an alternative discussion of this continuity criterion, cf. Tolstov (1976, chap. 9, particularly, sec. 9.1, p. 249) whom discusses it for a wave equation which satisfies Robin boundaries with one Dirichlet and one Neumann initial condition.

$$\dot{T}_n + \frac{D}{R}\lambda_n T_n = e^{st} f_n, \tag{6.4}$$

which is easily solved by multiplying by the integrating factor $e^{(D/R)\lambda_n t}$. Carrying out that multiplication, we will have $d/dt\,(e^{(D/R)\lambda_n t} T_n) = e^{(s+(D/R)\lambda_n)t} f_n$, so that setting $t = \tau$ and integrating for $t_0 < \tau < t$, we obtain

$$T_n(t) = e^{-(D/R)\lambda_n t} \int_{t_0}^{t} e^{(s+(D/R)\lambda_n)\tau} f_n(\tau)\,d\tau + e^{(D/R)\lambda_n(t_0-t)} T_n(t_0). \tag{6.5}$$

Since we have already provided a rationale for the assumption that the $T_n(t)$ are bounded, the series

$$w(x,t) = \sum_n T_n(t)\,\varphi_n(x) \tag{6.6}$$

automatically satisfies (6.2). Then, to satisfy the initial condition (4.12), we must have

$$\sum_n T_n(t_0)\,\varphi_n(x) = w(x,t_0), \tag{6.7}$$

from which we will isolate $T_n(t_0)$. Since the $\varphi_n(x)$ are continuously defined for $0 < x < \ell$, we can certainly multiply (6.7) by $\varphi_m(x)$, for any fixed $n = m$, to obtain $\sum_n T_n(t_0)\varphi_m\varphi_n = w(x,t_0)\varphi_m$. Since we have the assumption that the $T_n(t_0)$ will be at least bounded, it is reasonable to suppose that $\sum_n T_n(t_0)\varphi_m\varphi_n$ will converge uniformly to the continuous function $w(x,t_0)\varphi_m$ which permits the term-by-term integration[2]

$$\int_0^{\ell} w(x,t_0)\varphi_m(x)\,dx = \sum_n T_n(t_0) \int_0^{\ell} \varphi_m(x)\varphi_n(x)\,dx. \tag{6.8}$$

Referring to definition 5.14 for the inner product, the rhs of (6.8) may be written as $\sum_n T_n(t_0)(\varphi_m,\varphi_n)$. Since the eigenfunctions have the orthogonality $(\varphi_m,\varphi_n) = 0$ for all $n \neq m$, the only surviving term will be the one involving $n = m$, $T_m(t_0)(\varphi_m,\varphi_m)$, which is well defined for each fixed $m = 0,1,2,\ldots$. Thus, upon changing m to n, we have

$$T_n(t_0) = \frac{1}{(\varphi_n, \varphi_n)} \int_0^\ell w(x, t_0) \varphi_n(x) dx, \tag{6.9}$$

where we recall that the normalizing constants $(\varphi_n, \varphi_n)$ are positive real numbers.

---

[2]Proofs that the stated conditions are sufficient to admit a term-by-term integration can be found in (i) Widder (1989, chap. 9, sec. 6.2, theorem 13) in concise form, or in (ii) Knopp (1990, sec. 47, theorem **195**) as part of a larger expository discussion of term-by-term passage to the limit.

To complete the formula for the $T_n(t)$, and thereby a solution to the problem in w, we require the inhomogeneous term to have the expansion

$$F(x, \tau) = \sum_n f_n(\tau) \varphi_n(x). \tag{6.10}$$

We now multiply (6.10) by $\varphi_m$ to get $F\varphi_m = \sum_n f_n \varphi_m \varphi_n$, but now, since $F(x, t)$ is well defined and its temporal continuity will certainly pass to the time component of its decomposition $f_n(t)$, the series $\sum_n f_n \varphi_n \varphi_m$ must converge uniformly to the continuous function $F\varphi_m$ thus admitting the term-by-term integration

$$\int_0^\ell F(x, \tau) \varphi_m(x) dx = \sum_n f_n(\tau) \int_0^\ell \varphi_m(x) \varphi_n(x) dx. \tag{6.11}$$

We now refer to the discussion that intervenes between (6.8) and (6.9) for the fact that the rhs of (6.11) may be written as $\sum_n f_n (\varphi_m, \varphi_n)$. Then, by orthogonality, the only term that will remain is $f_m (\varphi_m, \varphi_m)$ whereupon, changing m to n, we obtain

$$f_n(\tau) = \frac{1}{(\varphi_n, \varphi_n)} \int_0^\ell F(x, \tau) \varphi_n(x) dx, \tag{6.12}$$

which completes the definition for $T_n(t)$ and for the problem in w.

We now begin our inquiry into the bound for $T_n(t)$ with the following lemma:

LEMMA 6.1 (SCHWARZ INEQUALITY): For two functions $h(\tau), k(\tau)$ real and continuous in the open interval $t_0 < \tau < t$ and bounded in the closed interval $t_0 \leq \tau \leq t$, we have the Schwarz inequality[3]

---

[3]For alternative discussions of the Schwarz inequality, cf. (i) Widder (1989, chap. 10, sec. 10.1, theorem 16) whose proof is based upon Cauchy's corresponding inequality for infinite series (ibid., chap. 9, sec. 8.1, theorem 16), and (ii) Tolstov (1976, sec. 2.4) whose proof is deduced from the observation that the polynomial cannot have distinct real zeros, and the inequality is shown to extend to the sum of any number of square integrable functions.

$$\left(\int_{t_0}^t h(\tau)\, k(\tau)\, d\tau\right)^2 \leq \int_{t_0}^t (h(\tau))^2\, d\tau \int_{t_0}^t (k(\tau))^2\, d\tau.$$

PROOF: Consider the expression $Q(\alpha) = (h + \alpha k, h + \alpha k)$ for $\alpha$ some real number. This will have an expansion as $(k,k)\alpha^2 + 2(h,k)\alpha + (h,h)$, a quadratic polynomial in $\alpha$. Completing the square, we obtain $[(2(k,k)\alpha + 2(h,k))^2 - \mathcal{D}]/(4(k,k))$, where $\mathcal{D}$ is the discriminant $\mathcal{D} = (2(h,k))^2 - 4(k,k)(h,h)$. Since we wish to determine a bound, and are therefore interested only in the non-negative values $|Q(\alpha)|$ which exist for $\mathcal{D} \leq 0$, we must have $(h,k)^2 \leq (h,h)(k,k)$.[4]

LEMMA 6.2 (BOUND FOR $T_n(t)$): Applying the Schwarz inequality to $T_n(t)$, we find:

1. For $n = 0$,

$$|T_0(t)| \leq \frac{e^{2st} - e^{(2/R)(Dr^2 t + \mu t_0)}}{2\mu/R(\varphi_0, \varphi_0)} \int_{t_0}^t \int_0^\ell (F(x,\tau))^2\, dx\, d\tau$$

$$+ \frac{e^{(2/R)(Dr^2 t + \mu t_0)}}{(\varphi_0, \varphi_0)} \int_0^\ell (e^{-rx}\phi(x) - H(x, t_0))^2\, dx.$$

(6.13)

2. For $n = 1, 2, 3, \ldots$,

$$|T_n(t)| \leq \frac{e^{2st} - e^{2(D/R)\lambda_n(t_0 - t) + 2st_0}}{2(s + (D/R)\lambda_n)(\varphi_n, \varphi_n)} \int_{t_0}^t \int_0^\ell (F(x,\tau))^2\, dx\, d\tau$$

$$+ \frac{e^{2(D/R)\lambda_n(t_0 - t) + 2st_0}}{(\varphi_n, \varphi_n)} \int_0^\ell (e^{-rx}\phi(x) - H(x, t_0))^2\, dx.$$

(6.14)

PROOF:

1. Setting $n = 0$ in (6.5), we define $T_0(t)$:
   (i) We begin with the first term, and referring to the Schwarz inequality, we set
   $h(\tau) = \left(e^{-(D/R)\lambda_0 t}\right) e^{(s + (D/R)\lambda_0)\tau}$, $k(\tau) = f_0(\tau)$ to obtain

---

[4] To verify that $(h,k)^2 \leq (h,h)(k,k)$ guarantees $|(h,k)|^2 \leq |(h,h)||(k,k)|$, begin with the non-negative expression $|Q(\alpha)| = |(h + \alpha k, h + \alpha k)|$ and, retracing each step in the proof, note that the absolute value operators drop out.

$$\frac{e^{2st} - e^{(2/R)(Dr^2 t + \mu t_0)}}{2\mu/R} \int_{t_0}^{t} (f_0(\tau))^2 \, d\tau$$

$$= \frac{e^{2st} - e^{(2/R)(Dr^2 t + \mu t_0)}}{2\mu/R} \int_{t_0}^{t} \left( \int_0^{\ell} \frac{\varphi_0}{(\varphi_0, \varphi_0)} F(x, \tau) \, dx \right)^2 d\tau,$$

where we have used the fact that $\lambda_0 = -r^2$, $s - Dr^2/R = \mu/R$, and $f_0(\tau)$ is from (6.12). We now apply the Schwarz inequality to the integral in parentheses. Setting $h(\tau) = \varphi_0/(\varphi_0, \varphi_0)$, we obtain the integrand $\varphi_0^2/(\varphi_0, \varphi_0)^2$ which can be rewritten as $(\varphi_0, \varphi_0)/(\varphi_0, \varphi_0)^2 = 1/(\varphi_0, \varphi_0)$ where, upon taking $k(\tau) = F(x, \tau)$, we arrive at the desired expression.

(ii) For the second term, we refer to (6.9) for $T_0(t_0)$, and taking $h(\tau) = \varphi_0/(\varphi_0, \varphi_0)$ with $k(\tau) = e^{(D/R)r^2(t-t_0)} w(x, t_0)$, we obtain

$$\frac{1}{(\varphi_0, \varphi_0)} \int_0^{\ell} \left( e^{(D/R)r^2(t-t_0)} w(x, t_0) \right)^2 dx$$

$$= \frac{e^{(2/R)(Dr^2 t + \mu t_0)}}{(\varphi_0, \varphi_0)} \int_0^{\ell} \left( e^{-rx} \phi(x) - H(x, t_0) \right)^2 dx,$$

where we have referred to (4.12) for $w(x, t_0)$.

2. Take $n = 1, 2, 3, \ldots$ in (6.5) to define $T_n(t)$:

(i) For the first term, allow $h(\tau) = \left( e^{-(D/R)\lambda_n t} \right) e^{(s + (D/R)\lambda_n)\tau}$ and perform the indicated integration. Then, taking $k(\tau) = f_n(\tau)$, we will have

$$\frac{e^{2st} - e^{2(D/R)\lambda_n(t_0 - t) + 2st_0}}{2(s + (D/R)\lambda_n)} \int_{t_0}^{t} (f_n(\tau))^2 \, d\tau$$

$$= \frac{e^{2st} - e^{2(D/R)\lambda_n(t_0 - t) + 2st_0}}{2(s + (D/R)\lambda_n)(\varphi_n, \varphi_n)} \int_{t_0}^{t} \int_0^{\ell} (F(x, \tau))^2 \, dx \, d\tau,$$

where the rhs follows immediately when $n = 1, 2, 3, \ldots$ is substituted for $n = 0$ in the discussion below the equation in part 1.i of this proof.

(ii) For the second term, referring to (6.9) for $T_n(t_0)$ and setting $h(\tau) = \varphi_n/(\varphi_n, \varphi_n)$ with $k(\tau) = e^{(D/R)\lambda_n(t_0-t)}w(x,t_0)$, we have

$$\frac{1}{(\varphi_n, \varphi_n)} \int_0^\ell \left(e^{(D/R)\lambda_n(t_0-t)}w(x,t_0)\right)^2 dx$$

$$= \frac{e^{2(D/R)\lambda_n(t_0-t)+2st_0}}{(\varphi_n, \varphi_n)} \int_0^\ell \left(e^{-rx}\phi(x) - H(x,t_0)\right)^2 dx,$$

which completes this proof.

COROLLARY 6.3 (BOUND FOR $T_n(t)$): The functions $T_n(t)$ are uniformly bounded in t, i.e., the sequence $\{T_n(t)\}$ is absolutely convergent for all t (however large).

PROOF: By (6.13),(6.14) $|T_n(t)|$ is bounded for any fixed $n,t$. Thus, it is only necessary to show that, for all t, the general term $|T_n(t)| \to 0$ with $n \to \infty$ and that, for all n, $|T_n(t)| < \infty$ at an arbitrarily large t, where it suffices to fix $t \geq 0$ and let $t_0 \to -\infty$.

(i) For $n=0$, we refer to (6.13), and noting that $|T_0(t)|$ does not depend upon n, we have

$$\lim_{t_0 \to -\infty} |T_0(t)| \leq \lim_{t_0 \to -\infty} \frac{e^{2st}}{2\mu/R(\varphi_0, \varphi_0)} \int_{t_0}^t \int_0^\ell (F(x,\tau))^2 \, dx \, d\tau < \infty.$$

(ii) For $n = 1, 2, 3, \ldots$, refer to (6.14) and observe that: (a) for any fixed t, $|T_n(t)| \to 0$ as $n \to \infty$; (b) for any fixed n, $|T_n(t)| \to 0$ as $t_0 \to -\infty$; and (c) $|T_n(t)| \to 0$ as $t_0 \to -\infty$, $n \to \infty$ simultaneously, noting that in each case $|T_n(t)|$ approaches zero in a continuous way.

By virtue of the fact that the $T_n(t)$ are bounded, the separated boundary conditions are continuously satisfied, and as a subtext of the arguments required to establish that bound, we found that the $T_n(t_0)$ are well defined for any n and any initial time $t_0$. The discussion of the auxiliary conditions will thus be complete if we can show that the initial condition (4.12) is continuously satisfied. For that we refer to (6.5) where, fixing $t_0$ and allowing t to approach $t_0$ from above, we see that the integral in the first term vanishes and the exponent in the second term approaches unity continuously. Therefore, the passage to the limit

$$\lim_{\substack{t \to t_0 \\ t > t_0}} T_n(t) \to T_n(t_0) \tag{6.15}$$

is termwise continuous.[5] This will be sufficient to ensure that the initial condition is continuously satisfied if we can establish the convergence required to justify the passage from (6.7) to (6.8).

THEOREM 6.4 (CONVERGENCE OF $\sum_n T_n \varphi_n$): The infinite series $\sum_n T_n \varphi_n$ is uniformly and unconditionally convergent.[6]

PROOF: Begin by noting that $T_n(t)\varphi_n(x) = w_n(x,t)$ is continuously defined for $\{x : 0 \leq x \leq \ell\}$, $\{t : -\infty < t_0 \leq t, t > 0\}$. Then by Weierstrass,[7] it is sufficient to show that $|w_n(x,t)| \leq M_n$ where $M_n \geq 0$ is a real number belonging to the convergent series $\sum_n M_n$.

(i) For $n = 0$, we require $|w_0(x,t)| = |T_0(t)| e^{rx} \leq M_0$. By lemma 6.2.1 and part (i) of the proof of corollary 6.3, $|T_0(t)|$ will be a number for any fixed t, from which we obtain the bound $M_0 = |T_0(t)| e^{r\ell}$.

(ii) For $n = 1, 2, 3, \ldots$, the requirement that $|w_n(x,t)| = |T_n(t)||\mathrm{Cos}(x\sqrt{\lambda_n}) + (r/\sqrt{\lambda_n}) \mathrm{Sin}(x\sqrt{\lambda_n})| \leq M_n$ will be satisfied if $\sum_n M_n = \sum_n |T_n(t)| (1 + r/\sqrt{\lambda_n}) < \infty$. To confirm that this series converges, we refer to (6.14) where we fix $t > 0$ and note that $t_0 - t$ is consequently a negative numeric exponent. We then carry out the multiplication by $1 + r/\sqrt{\lambda_n}$ and recall that $\sqrt{\lambda_n} = n\pi/\ell$ and that $(\varphi_n, \varphi_n) = (r^2 + \lambda_n)\ell/(2\lambda_n)$. When the resulting expression is fully simplified, the first term in (6.14) converges by comparison to the p series $\alpha/n^p$ where $p = 2$ with $\alpha = e^{2st}, e^{2st_0}$ the appropriate constant. The second term in (6.14) converges with its d'Alembert ratio approaching zero exponentially fast with increasing n. To complete the proof, we note that, by virtue of part (ii) of the proof of corollary 6.3, the convergence we have established will hold for all t.

---

[5] Note that the limit in (6.15) is the product of a careful formulation which in no way implies that we may let $t \to t_0$ from $t > t_0$ and get $w(x,t) \to w(x,t_0)$. In fact, because the latter limit usually fails to yield anything meaningful, the engineering literature refers to diffusive processes as irreversible. In the mathematics literature, this irreversibility is referred to as the backward-time problem, which is shown to yield a solution that is generally non-unique. For a discussion of the backward-time problem in (i) a finite interval, with a sketch of the proof of non-uniqueness, see Zachmanaglou and Thoe (1986, chap. 9, sec. 2, especially the discussion on p. 339 which follows from theorem 2.1), (ii) an infinite interval, see Berg and McGregor [1966, sec. 9.6, notably the discussion on p. 252 supporting their eqs. (9.6.16)-(9.6.18)].

[6] Unconditional convergence has the usual meaning that the series will converge no matter how the terms are arranged, as guaranteed when the terms converge in absolute value (for a discussion of the many theorems which accrue to unconditionally convergent series, see Knopp 1990, sec. 16).

[7] This is referred to as Weierstrass's M-test. For two different proofs, both carried out in the context of a larger discussion of uniform convergence, see (i) Widder (1989, chap. 9, sec. 5, particularly sec. 5.2), or (ii) Knopp (1990, chap. 11, especially sec. 48, theorem **197**).

## 6.2. SOLUTION TO THE ORIGINAL SYSTEM

In this section, we furnish the final elements for a solution to the system specified in section 1.2. Since the intervening discussion has provided a well-defined representation for w, a well-defined expression for C follows immediately from the inverse linear transformation

$$C(x,t) = \left(w(x,t) + e^{st}H(x,t)\right)e^{rx-st} \quad (6.16)$$

of chapter 4, where H is the known function (4.9) and we have referred to (4.8) for $u(x,t) = w(x,t) + e^{st}H(x,t)$ and to (4.1),(4.3) for $C(x,t) = u(x,t)e^{rx-st}$.

To complete the representation, it only remains for us to provide a large-t solution and a formula for the exit concentration. For the large-t solution, we refer to appendix-B (3.5) which is derived from appendix-B (2.5) which is identical to (6.16). To provide an expression for the exit concentration defined in section 1.2, we refer to appendix B, section 4, which develops the transformations required to furnish a canonical diffusion equation thus admitting a real-valued solution.

**CHAPTER 7. CONCLUDING DISCUSSION**

7.1. ON THE WELL-POSEDNESS OF THE SYSTEM

The question of existence has been addressed constructively although we have yet to establish uniqueness and stability:

To establish uniqueness, recall from section 4.3 that our separated solution generates the linear system SLR and from the closing paragraph of chapter 5 that any subset consisting of N-many solutions of SLR will have the Hamel basis $\{\varphi_n\}_0^N$. Thus, taking $t_i$ some fixed point in the interval $t_0 < t_i \leq t$, we may define the constant $T_n(t_i) = \alpha_{n,t_i}$ so that $\{\alpha_{n,t_i}\varphi_n\}_0^N$ will be a linear combination of solutions to the system composed by the ODE (6.3) and boundary conditions (6.2).[1] Then, since $t_i$ is arbitrary, we may let $t_i = t$ so that $\alpha_{n,t} = T_n(t)$, and since we may choose N however we please with $\sum_{N+1}^\infty T_n \varphi_n \to 0$ uniformly, the Hamel basis $\{T_n\varphi_n\}_0^N$ will furnish a uniform approximation to w. The well-posedness of the system articulated in section 1.2 then follows immediately from the linearity of the inverse transformations discussed in section 6.2.

To establish stability means to show that a small change in the data will cause only a small change in the solution, which follows easily from our assumptions on the physical problem. Since we require well behaved initial and boundary concentrations and allow only smooth changes in the physical data, our solution converges uniformly to its defining function. Since our solution is analytic in all of its variables, any permissible change in $\phi, g$ is obviously smoothed into the solution, thus verifying physical intuition for a diffusive process.

7.2. ON THE NECESSITY OF ROBIN BOUNDARY CONDITIONS

THEOREM 7.1: Robin boundaries are a necessary and sufficient condition of any solution to the CDE.

PROOF: (i) Sufficiency (mass conservation): Consider the overall mass balance illustrated in figure 7.1.

---

[1] Since we have established that $T_n(t)$ is uniformly bounded for all t, there will always be some (real) number $K \geq 0$ such that $|\alpha_{n,t_i}| \leq K$ with the consequence that $w_{n,t_i}(x) = \alpha_{n,t_i}\varphi_n(x)$ will be a simple linear combination of the $\varphi_n$. Then returning to chapter 5, we determine an operator for (6.3) similar to S and verify its linearity with lemma 5.2 and corollary 5.22. Now, note that $w_{n,t_i}(x)$ will satisfy the boundary conditions (6.2) by satisfying the linear operators $\alpha_{n,t_i}B_1(\varphi_n)$, $\alpha_{n,t_i}B_2(\varphi_n)$. The superposition $w_{N,t_i}(x) = \sum_0^N \alpha_{n,t_i}\varphi_n(x)$ thus solves the differential equation (6.3) with end conditions (6.2) for any given $t = t_i$.

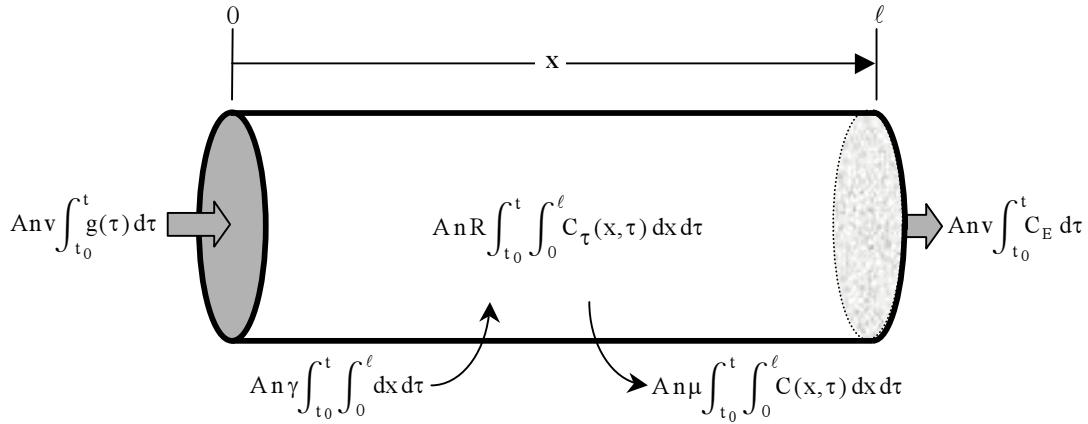

FIGURE 7.1. OVERALL MASS BALANCE.

The resulting conservation statement is

$$R \int_{t_0}^{t} \int_{0}^{\ell} C_\tau(x,\tau) \, dx \, d\tau = v \int_{t_0}^{t} \left(g(\tau) - C_E\right) d\tau + \int_{t_0}^{t} \int_{0}^{\ell} \left(\gamma - \mu C(x,\tau)\right) dx \, d\tau. \qquad (7.1)$$

Since, we have shown that $C_\tau$ is continuous in $(t_0, t)$, $(0, \ell)$, we may interchange the order of integration on the lhs and integrate over time to obtain[2]

$$R \int_{0}^{\ell} \left(C(x,t) - \phi(x)\right) dx = v \int_{t_0}^{t} \left(g(\tau) - C_E\right) d\tau - \mu \int_{t_0}^{t} \int_{0}^{\ell} C(x,\tau) \, dx \, d\tau + \gamma \ell (t - t_0),$$

where we have used the initial condition $C(x, t_0) = \phi(x)$. We now derive a mathematical statement by integrating (1.1) over the length of the column to get

$$R \int_{0}^{\ell} C_t(x,t) dx = D \int_{0}^{\ell} C_{xx}(x,t) dx - v \int_{0}^{\ell} C_x(x,t) dx + \int_{0}^{\ell} \left(\gamma - \mu C(x,t)\right) dx,$$

and carrying out the integration for the first two terms on the rhs will yield

---

[2]For a discussion of the conditions under which integration is commutative, see Widder's (1989, chap. 6, sec. 3, exercise 11, and chap. 6, sec. 5.1) discussion of Dirichlet's formula for Riemann integrals.

$$R \int_0^\ell C_t(x,t)\,dx \\ = DC_x(\ell,t) - vC(\ell,t) - \left(DC_x(0,t) - vC(0,t)\right) + \int_0^\ell (\gamma - \mu C(x,t))\,dx. \tag{7.2}$$

Upon returning to the conservation equation (7.1) and evaluating its time derivative, we find that[3]

$$R \int_0^\ell C_t(x,t)\,dx = vg(t) - vC_E + \int_0^\ell (\gamma - \mu C(x,t))\,dx. \tag{7.3}$$

If we now equate the physical statement (7.3) to the mathematical statement (7.2), we will have $vC(0,t) - DC_x(0,t) - \left(vC(\ell,t) - DC_x(\ell,t)\right) = vg(t) - vC_E$, and upon rearranging, we obtain

$$vC(0,t) - DC_x(0,t) - vg(t) = vC(\ell,t) - DC_x(\ell,t) - vC_E. \tag{7.4}$$

Note that the lhs is our entrance condition and the rhs the exit condition so that (7.4) is the identity $0 \equiv 0$. Thus, solutions conserve mass when they satisfy Robin boundaries, which verifies their sufficiency.

(ii) Necessity (internal-consistency): Let $C(x,t)$ be a solution to the CDE that correctly describes the inlet and exit concentrations (and therefore conserves mass). Now let $C(x,t)$ satisfy a continuous-concentration hypothesis at either boundary. Since we have shown that $C(x,t)$ must be continuously defined in $0 \leq x \leq \ell$, we can fix x at either boundary and let the differential element $\Delta x$ in (3.2) approach the boundary surface from the inside. As $\Delta x$ becomes infinitesimal, the derivative must vanish to satisfy our hypothesis, but then $(C_x(x+\Delta x,t) - C_x(x,t))/\Delta x$ must also vanish so that $C(x,t)$ must not be a solution to the CDE, and this completes the theorem.[4]

7.3. CONCLUSION

Our conclusion can now be summarized in a short statement: By section 7.1 the solution we have provided is unique, and by section 7.2 we have solved the only system that is both physically and mathematically correct. It is therefore impossible to derive another solution to the stated CDE that differs in any material way.

---

[3] See Widder (1989, chap. 10, sec. 7.1, theorem 8) for a discussion of differentiation of an integral with respect to its upper limit, which is easily permitted for our well-behaved solutions.

[4] Theorem 7.1 has a corollary for the historic boundary conditions (see app. C).

**APPENDIX A. LIST OF SYMBOLS (UNITS)**

- $\alpha$    dispersivity (L)
- $\mu$    first-order decay rate $(T^{-1})$
- $\phi$    initial concentration $(ML^{-3})$
- $\ell$    fixed length of the flow domain (L)
- $\rho$    porous-medium bulk density $(ML^{-3})$
- $\gamma$    constant production-rate $(ML^{-3}T^{-1})$
- A    area $(L^2)$
- C    solute concentration $(ML^{-3})$
- D    hydrodynamic dispersion $(L^2T^{-1})$
- g    solute input concentration $(ML^{-3})$
- k    linear distribution coefficient $(L^3M^{-1})$
- n    porous-medium intrinsic porosity (unitless)
- q    specific discharge-rate $(LT^{-1})$
- R    retardation factor (unitless)
- t    time (T)
- v    average interstitial fluid velocity $(LT^{-1})$
- x    distance (L)

SUBSCRIPTS

- d    diffusion (Fickian)
- E    exit
- F    flux
- L    liquid
- S    sorbed
- t    partial derivative in t
- x    partial derivative in x
- xx    second partial derivative in x

**Appendix B. The Convection-Diffusion Equation for a Finite Domain with Time Varying Boundaries (preprint version)**

# The convection-diffusion equation for a finite domain with time varying boundaries[1,2,3]

W. J. GOLZ
Department of Civil and Environmental Engineering, Louisiana State University
Baton Rouge, Louisiana 70803, USA

J. R. DORROH
Department of Mathematics, Louisiana State University
Baton Rouge, Louisiana 70803, USA

**Abstract**—A solution is developed for a convection-diffusion equation describing chemical transport with sorption, decay, and production. The problem is formulated in a finite domain where the appropriate conservation law yields Robin conditions at the ends. When the input concentration is arbitrary, the problem is underdetermined because of an unknown exit concentration. We resolve this by defining the exit concentration as a solution to a similar diffusion equation which satisfies a Dirichlet condition at the left end of the half line. This problem does not appear to have been solved in the literature, and the resulting representation should be useful for problems of practical interest.

Authors of previous works on problems of this type have eliminated the unknown exit concentration by assuming a continuous concentration at the outflow boundary. This yields a well-posed problem by forcing a homogeneous Neumann exit, widely known as the Danckwerts [1] condition. We provide a solution to the Neumann problem and use it to produce an estimate which demonstrates that the Danckwerts condition implies a zero concentration at the outflow boundary, even for a long flow domain and a large time.

**Keywords**—Chemical transport, Danckwerts conditions, Error estimate, Robin conditions.

## 1. INTRODUCTION

CLASSICAL SOLUTIONS of the diffusion equation have been catalogued for many of the important problems of heat transfer [2]. The diffusion equation has also been widely employed as a model for chemical reaction processes, and this usually entails the inclusion of lower-order terms that describe convection and reaction. The complications introduced by those additional terms often call for inventive techniques which yield novel and useful representations (e.g., see [3]).

Descriptions of chemical processes and contaminant transport have motivated a large volume of work on the one-dimensional convection-diffusion equation (CDE). Numerical methods have provided solutions to problems satisfying a fairly wide range of conditions (e.g. see [3-5]). However, analytical solutions continue to be highly valued for their inherent simplicity, their capacity to convey qualitative information about the physical problem, and as a verification for numerical models. Many of the existing analytical solutions have been compiled in two large compendiums [6,7].

---

[1] ERRATUM: Equation (2.4) should be corrected to read $vC(\ell,t) - DC_x(\ell,t) = vC_E$ (this erratum was corrected in the published version of this paper, for which the citation is provided in footnote 2).

[2] FOR CITATION, see the final published version of this paper: Golz, W. J. and J. R. Dorroh. 2001. The Convection-diffusion equation for a finite domain with time varying boundaries. *Applied Mathematics Letters* 14 : 983-988 (received by *AML* September 2000; accepted by *AML* October 2000).

[3] ERRATA: (i) In equation (2.7), the term $-H_t$ should be corrected to read $-(sH + H_t)$; (ii) In equation (2.9), the exponential factor $e^{-\ell x}$ should be corrected to read $e^{-r\ell}$ (these errata exist in this preprint and in the published version of this paper).

The literature on the mathematics of chemical transport demonstrates that problems posed with a Neumann condition fail to satisfy the relevant conservation equation for concentrations in the interior of the domain [8-10]. The inconsistency with the physical problem arises because the conservation of mass requires Robin conditions at the ends, which leads to a system that is underdetermined due to an unknown exit concentration [1,10]. Danckwerts eliminated the unknown exit concentration by assuming a continuous concentration at the outflow boundary, thus forcing a homogeneous Neumann condition [1,10]. The Danckwerts assumption yields a well-posed problem, but it introduces error into predicted concentrations [6,10]. The magnitude of the error is inversely proportional to time and Peclet number (i.e., the product of fluid velocity and flow length when divided by the diffusion constant, see [6,10]). Hydrogeological problems characteristically exhibit a low Peclet number, and observations at early times are often important. Thus, the error can be significant, especially near an outflow boundary [6,10].

This paper was motivated by a class of physically important problems from contaminant hydrogeology. These problems can be characterized by a diffusion equation with additional terms for drift, decay, and production, and the input may be any bounded continuously differentiable function. The resulting system is composed of an inhomogeneous CDE with time varying Robin conditions. Mass is conserved when the exit concentration is defined as the solution to a CDE with a Dirichlet condition at the left end of the half line. This approach yields the proper formulation in finite geometry thus admitting a representation in an eigenfunction expansion with closed-form eigenvalues.

## 2. FORMULATION AND TRANSFORMATION

The problems that will be considered in this paper take the form

$$RC_t = DC_{xx} - vC_x - \mu C + \gamma, \qquad 0 < x < \ell, \quad t \in \mathbb{R} \tag{2.1}$$

where $C = C(x,t)$ denotes a concentration, with x the longitudinal distance and t the time. The constants R, D, v, μ, and γ describe linear sorption, diffusion, longitudinal fluid velocity, decay, and production, respectively. Equation (2.1) will satisfy the auxiliary conditions

$$C(x, t_0) = \phi(x), \qquad 0 < x < \ell \tag{2.2}$$

$$vC(0,t) - DC_x(0,t) = vg(t), \qquad t > t_0 \tag{2.3}$$

$$vC(\ell, t) - DC_x(\ell, t) = C_E, \qquad t > t_0 \tag{2.4}$$

where $\phi(x), g(t)$ are bounded continuously differentiable functions, and $C_E$ is the exit concentration which, for now, we regard as an experimentally measured quantity.

To transform the CDE into standard form for an inhomogeneous diffusion equation, and to provide homogeneous boundary conditions, we introduce the change of dependent variable

$$C(x,t;r,s) = \left(w(x,t) + e^{st} H(x,t)\right) e^{rx-st} \tag{2.5}$$

If we choose the parameters as

$$r = \frac{v}{2D}, \qquad s = \frac{1}{R}\left(\frac{v^2}{4D} + \mu\right), \qquad D, R \neq 0 \tag{2.6}$$

then (2.1) may be written as the inhomogeneous diffusion equation

$$w_t = \frac{D}{R} w_{xx} + e^{st} F(x,t), \qquad F(x,t) = \frac{\gamma}{R} e^{-rx} - H_t + \frac{D}{R} H_{xx} \tag{2.7}$$

and the initial condition (2.2) will become

$$w(x, t_0) = e^{st_0}\left(e^{-rx}\phi(x) - H(x, t_0)\right) \tag{2.8}$$

If we then define

$$H(x, t) = \left(1 + \cos\frac{\pi x}{\ell}\right)g(t) + \left(1 - \cos\frac{\pi x}{\ell}\right)e^{-\ell x}C_E \tag{2.9}$$

the boundary conditions (2.3), (2.4) will be homogeneous:

$$w_x(0, t) - r\,w(0, t) = 0 \tag{2.10}$$

$$w_x(\ell, t) - r\,w(\ell, t) = 0 \tag{2.11}$$

## 3. EIGENFUNCTION EXPANSION

To obtain our representation, we only require the separated solution $w(x,t) = \varphi(x)T(t)$ to be justified when $F(x,t) = 0$. This leads to a regular Sturm-Liouville problem with the positive eigenvalues and related eigenfunctions

$$\lambda_n = \frac{n^2 \pi^2}{\ell^2}, \qquad \varphi_n(x) = \cos\left(x\sqrt{\lambda_n}\right) + \frac{r}{\sqrt{\lambda_n}}\sin\left(x\sqrt{\lambda_n}\right), \qquad n = 1, 2, 3, \cdots \tag{3.1}$$

where each $\lambda_n$ defines a single-dimensional eigenspace. There is also one negative eigenvalue with a single associated eigenfunction:

$$\lambda_n = -r^2, \qquad \varphi_n(x) = e^{rx}, \qquad n = 0 \tag{3.2}$$

The negative eigenvalue arises from the transformations, so its related eigenfunction does not entail the usual meaning of a source term which would allow the solution to become arbitrarily large as t increased. In fact, we will exhibit a large-t solution that is bounded.

The normalized eigenfunctions of a Sturm-Liouville system clearly compose an orthonormal basis (e.g., see [11], sec. 7.5), from which it follows that (2.7) will have a solution of the form

$$w(x, t) = \sum_{n=M}^{\infty} \varphi_n(x)\, T_n(t), \qquad M = 0 \tag{3.3}$$

where

$$T_n(t) = \frac{e^{-\frac{D}{R}\lambda_n t}}{\int_0^\ell \varphi_n^2(x)\,dx}\left(\int_{t_0}^t e^{\left(s + \frac{D}{R}\lambda_n\right)\tau}\int_0^\ell F(x,\tau)\varphi_n(x)\,dx\,d\tau \right.$$

$$\left. + e^{\left(s + \frac{D}{R}\lambda_n\right)t_0}\int_0^\ell \left(e^{-rx}\phi(x) - H(x, t_0)\right)\varphi_n(x)\,dx \right) \tag{3.4}$$

To obtain (3.4), the orthogonality of the eigenfunctions and the symmetry of the boundary conditions on w and φ are required to develop an ordinary differential equation in T that satisfies the initial condition in w (for a discussion of this technique, see [12], chap. 9).

We now wish to demonstrate that the representation of C(x,t) will be bounded for all t. To establish this, we fix t in (3.4) and then note that the second term is bounded and will therefore vanish as $t_0 \to -\infty$. Then referring to (2.5), a large-t solution to (2.1)-(2.4) will be

$$C_{bv}(x,t) = \sum_{n=M}^{\infty} \frac{\varphi_n(x) e^{rx - \left(s + \frac{D}{R}\lambda_n\right)t}}{\int_0^\ell \varphi_n^2(x)\,dx} \int_{-\infty}^t e^{\left(s + \frac{D}{R}\lambda_n\right)\tau} \int_0^\ell F(x,\tau)\varphi_n(x)\,dx\,d\tau + e^{rx} H(x,t) \quad (3.5)$$

which is obviously nonincreasing for t. Note that (3.5) is equivalent to a solution of the pure boundary value problem (2.1), (2.3), (2.4) where it is valid for all $t \in \mathbb{R}$ (for a similar result involving a problem on the half line, see [3]).

## 4. EXIT CONCENTRATION

In the event that the exit concentration is not empirically observed, we may determine it. First, set

$$vC_F(x,t) = vC(x,t) - DC_x(x,t), \qquad 0 \le x < \infty, \quad t \in \mathbb{R} \quad (4.1)$$

It can be verified by direct substitution that $C_F$ satisfies a differential equation identical to (2.1) (e.g., see [8,9]). It then follows that $C_E = C_F(\ell, t)$, where $C_F(x, t)$ is a bounded solution to (2.1) with the auxiliary conditions

$$vC_F(x,t_0) = v\phi(x) - D\phi_x(x), \qquad 0 < x < \infty \quad (4.2)$$

$$C_F(0,t) = g(t), \qquad t > t_0 \quad (4.3)$$

The change of variable

$$C_F(x,t) = u(x,t) e^{rx - st} + \frac{\gamma}{\mu} \quad (4.4)$$

will, with reference to (2.1), yield the canonical form of the diffusion equation

$$u_t = \frac{D}{R} u_{xx} \quad (4.5)$$

and (4.2), (4.3) will become

$$u(x,t_0) = \Phi(x) \quad (4.6)$$

$$u(0,t) = G(t) \quad (4.7)$$

where

$$\Phi(x) = \left(\phi(x) - \frac{D}{v}\phi_x(x) - \frac{\gamma}{\mu}\right) e^{-rx + st_0}, \qquad G(t) = \left(g(t) - \frac{\gamma}{\mu}\right) e^{st} \quad (4.8)$$

A solution to (4.5)-(4.7) is provided by the familiar integral representation

$$u(x,t) = \int_0^\infty K(x-\zeta,(D/R)(t-t_0)) - K(x+\zeta,(D/R)(t-t_0))\Phi(\zeta)\,d\zeta \\ - \frac{2D}{R}\int_{t_0}^t K_x(x,(D/R)(t-\tau))G(\tau)\,d\tau \quad (4.9)$$

where the first term may be obtained from an odd extension of $\Phi$ and the last term from the method of Duhamel (e.g., see [2], chap. 4). The kernel and its derivative are then

$$K(x,t) = \frac{1}{\sqrt{4\pi t}} e^{-x^2/4t}, \qquad K_x(x,t) = -\frac{x}{2t} K(x,t) \quad (4.10)$$

and may be obtained in a straightforward manner from a Fourier-transform solution of the fundamental initial-value problem (e.g., see [2], chap. 3) or from a Green's function (e.g., see [3]).

## 5. AN ESTIMATE OF THE DANCKWERTS ERROR

In this section, we provide a solution to the Danckwerts problem (i.e., the Robin exit condition is replaced with a Neumann condition). The resulting representation is then used to produce an estimate of the Danckwerts error.

First we define a new concentration $C_D = C_D(x,t)$ which will satisfy (2.1)-(2.3), but (2.4) will be replaced by the Neumann condition

$$C_{Dx}(\ell, t) = 0, \qquad t > t_0 \tag{5.1}$$

To see that (5.1) is not generally valid, take $C_{Dt} = \gamma = 0 \ni C_D = C_D(x)$ in (2.1) which yields the problem considered by Danckwerts [1]. Now note that $C_{Dx}(\ell) = 0 \Rightarrow C_D(\ell) = 0$.

The change of variable (2.5) and parameters (2.6) remain the same while the CDE (2.7) and the initial condition (2.8) must respect the new definition

$$H_D(x,t) = \left(1 + \cos\frac{\pi x}{\ell}\right) g(t) \tag{5.2}$$

which replaces (2.9). The entrance condition (2.10) for $w_D = w_D(x,t)$ will remain unchanged, but the exit condition (2.11) will now be replaced by

$$w_{Dx}(\ell, t) + r w_D(\ell, t) = 0 \tag{5.3}$$

The separated solution $w_D(x,t) = \varphi_D(x) T_D(t)$ for the homogeneous problem in $w_D$ will now yield a Sturm-Liouville problem with strictly positive eigenvalues. The eigenvalues will no longer be given in the closed form of (3.1) but will instead arise from the sequential intersections $z_1 = z_2$, where

$$z_1 = \tan\left(\ell \sqrt{\lambda_D}\right), \qquad z_2 = \frac{2r\sqrt{\lambda_D}}{\lambda_D - r^2} \tag{5.4}$$

Information about how the $\lambda_D$ approach their asymptotic value is useful for evaluating the error estimate. First, observe that

$$\frac{n^2 \pi^2}{\ell^2} < \lambda_{Dn} < \frac{(n+1)^2 \pi^2}{\ell^2} \tag{5.5}$$

so for each n, $\lambda_n < \lambda_{Dn}$. Then, the $\lambda_{Dn}$ will have the asymptotic value

$$\lim_{\sqrt{\lambda_D} \to \infty} z_2 = \lim_{\sqrt{\lambda_D} \to \infty} \frac{2r\sqrt{\lambda_D}}{\lambda_D - r^2} = \lim_{\sqrt{\lambda_D} \to \infty} \frac{r}{\sqrt{\lambda_D}} = 0 \Rightarrow \lim_{n \to \infty} \lambda_{Dn} = \frac{n^2 \pi^2}{\ell^2} \tag{5.6}$$

The eigenfunctions $\varphi_{Dn}(x)$ will be defined as in (3.1) but with their eigenvalues now given by (5.4). If we let the summation begin at $M = 1$, a solution for $w_D(x,t)$ will be provided by (3.3) where the $T_{Dn}(t)$ come from (3.4), and a large-t solution will be given by $C_{Dbv}(x,t)$ as in (3.5).

Since C and $C_D$ differ most near $x = \ell$, our definition of the error is the natural one

$$E_D(t) = |C(\ell, t) - C_D(\ell, t)| \tag{5.7}$$

If we allow that (5.7) is, by hypothesis, nonincreasing with t, it then follows that

$$E_D(t) \geq |C_{bv}(L, t) - C_{Dbv}(L, t)| \geq |\gamma/\mu| \tag{5.8}$$

The verification of (5.8) is left to the reader, but it is easy to see that $|\lambda/\mu|$ will be obtained from $|C_{bv}(\ell, t) - C_{Dbv}(\ell, t)|$ at a sufficiently large constant length L since $\lambda_n, \lambda_{Dn} \to 0$ for the first few eigenvalues and similar terms in the large-t solutions $C_{bv}(\ell, t), C_{Dbv}(\ell, t)$ approach one another quite rapidly as n increases.

We now refer to (2.1), where it is clear that $|\lambda/\mu|$ is the amount that C(x,t) will differ from zero at a very large distance. Thus, the Neumann condition implies a zero concentration at the outflow boundary, just as in the original problem considered by Danckwerts [1].

**APPENDIX C. COROLLARY OF THEOREM 7.1: HISTORIC BOUNDARY CONDITIONS**

Theorem 7.1 has a corollary for the historic boundary conditions that employ continuous-concentration hypotheses. Whether continuity assumptions are rendered as a Dirichlet or a homogeneous-Neumann condition, they satisfy mass-conservation condition 7.1.i by matching their solutions to exterior concentrations. To achieve that, their concentrations must sacrifice internal consistency with the CDE and will therefore fail to satisfy 7.1.ii.

Consider the specific case of the Danckwerts' problem discussed in section 5 of appendix B. Begin by referring to appendix-B (5.1), noting that the preliminary transformations for the Danckwerts' problem yield expressions identical to those in section 4.1 of the dissertation. However, to transform the Danckwerts' boundaries to a homogeneous form, we use appendix-B (5.2) in place of dissertation (4.9). We then carry out the transformations of the Danckwerts' problem up to dissertation (4.11) where we find that $F_2 \equiv 0$. Thus, the Neumann condition furnishes a solution to the CDE only when the exit concentration is identically zero. Because that can be generally true only for steady-state problems, and then only when no unreacted mass reaches the exit, Danckwerts' solutions contain the error described by estimates (5.7),(5.8) of appendix B.

To see that the continuous-concentration hypotheses employed by WW yield an equivalent error, simply return to the discussion in the literature review and note that a WW solution will coincide with that of Danckwerts within the reactor.